\documentclass[a4paper,11pt]{article}

\usepackage[pdftex]{hyperref}
\usepackage{stmaryrd}
\usepackage[normalem]{ulem}
\usepackage{graphicx,xcolor}
\usepackage{amsfonts}
\usepackage{amsmath,amssymb,mathrsfs,amsthm}
\usepackage[all]{xy} 
\usepackage{color}
\usepackage{times}
\usepackage{enumerate,vmargin}
\usepackage{verbatim}
\usepackage{paralist}
\usepackage[labelfont=bf,size=small,font=it]{caption}
\usepackage[numbers,sort&compress]{natbib}
\usepackage{tikz}
\usetikzlibrary{shapes,arrows} 

\setpapersize{A4}
\setmarginsrb{2.5cm}{2cm}{2.5cm}{2cm}{.4cm}{4mm}{0cm}{7.5mm}
 
\hypersetup{
    unicode      = false,     
    pdftoolbar   = true,      
    pdfmenubar   = true,      
    pdffitwindow = true,      
    pdfnewwindow = true,      
    colorlinks   = true,      
    linkcolor    = blue,      
    citecolor    = blue,      
    filecolor    = blue,      
    urlcolor     = blue       
}

\theoremstyle{plain}
\newtheorem{lem}{Lemma}
\newtheorem{thm}[lem]{Theorem}
\newtheorem{prop}[lem]{Proposition}

\newtheorem*{thm*}{Theorem 1*}

\theoremstyle{definition}

\theoremstyle{remark}

\newcommand{\cR}{\mathcal R}

\newcommand{\bP}{\mathbf P}

\newcommand{\bp}{\mathbf p}

\newcommand{\sL}{\mathscr L}

\newcommand{\llb}{\llbracket}
\newcommand{\rrb}{\rrbracket}

\newcommand{\Br}{\operatorname{Br}}

\newcommand{\Skel}{\operatorname{Skel}}
\newcommand{\Card}{\operatorname{Card}}

\newcommand{\cT}{\mathcal T}

\newcommand{\cC}{\mathcal C}
\newcommand{\cP}{\mathcal P}
\newcommand{\cS}{\mathcal S}

\newcommand{\cL}{\mathcal L}

\newcommand{\bth}{\boldsymbol{\theta}}

\newcommand{\R}{\mathbb R}
\newcommand{\N}{\mathbb N}

\newcommand{\eqd}{\stackrel{d}=}

\usepackage{xcolor}
\usepackage[many]{tcolorbox}
\tcbuselibrary{breakable}

\graphicspath{{.}{figures/}}

\title{\textsc{Pruning, Cut trees, and the Reconstruction Problem}}

\author{Nicolas \textsc{Broutin}
\thanks{LPSM, Sorbonne Universit\'e, Campus Pierre et Marie Curie, 
4 place Jussieu,  
75252 Paris Cedex 05,  
France. Email: nicolas.broutin@sorbonne-universite.fr} 
\thanks{Institut Universitaire de France (IUF)}
\and Hui \textsc{He}
\thanks{School of Mathematical Sciences, Beijing Normal University, Beijing 100875,
China.
Email: hehui@bnu.edu.cn}
\and Minmin \textsc{Wang}
\thanks{Department of Mathematics, School of Mathematical and Physical Sciences, University of Sussex, 
Brighton, BN1 9QH,  United Kingdom.  
Email: minmin.wang@sussex.ac.uk }
}

\date{\today}

\begin{document}

\maketitle
\begin{abstract}
We consider a pruning of  the inhomogeneous continuum random trees, as well as the cut trees that encode the genealogies of the fragmentations that come with the pruning. We propose a new approach to the reconstruction problem, which has been treated for the Brownian CRT in \cite{BrWa} [\emph{Electron.\ J.\ Probab.\ vol.\ 22, 2017}] and for the stable trees in \cite{ADG} [\emph{Ann. IHP B, vol 55, 2019}]. Our approach does not rely upon self-similarity and can potentially apply to general L\'evy trees as well. 




\end{abstract}



\section{Introduction}

In this work, we consider pruning processes of random tree-like metric spaces, which arise in the scaling limits of prunings of finite discrete trees. Extending the notion of tree as a loop-free connected graph, we say that a complete metric space $(\cT,d)$ is a {\it real tree} if for any two points $\sigma,\sigma'\in \cT$ there is a unique arc with end points $\sigma$ and $\sigma'$ denoted by $\llb \sigma,\sigma'\rrb$, and that arc is isometric to the interval $[0, d(\sigma, \sigma')]$. 
The reader is encouraged to refer to \cite{Evans} for a comprehensive study on real trees. 
We are interested in random real trees, or more precisely a collection of fractal-like random real trees called {\it continuum random trees} (CRT in short). 

Let us take a special case in order to illustrate our general motivation: the celebrated Brownian continuum random tree of Aldous \cite{Aldous1991b}. 
Aldous and Pitman \cite{Al00} have defined a fragmentation process related to the standard additive coalescent. Informally, this fragmentation process is engineered by logging the Brownian CRT at uniform locations on its skeleton, resulting in smaller and smaller connected components as time goes on. 
There is a natural notion of genealogy associated to this fragmentation, which can be represented as a separate continuum random tree, called the \emph{cut tree}. It turns out that this cut tree is also distributed like a Brownian continuum random tree \cite{BeMi13}. This distributional identity between the CRT and its cut tree in fact holds for more general families of random trees, including the stable L\'evy trees \cite{Di15}, the general L\'evy tree under the excursion measures \cite{AbDe13}, as well as the inhomogeneous continuum random trees \cite{BW}. 

Roughly speaking, the reconstruction problem asks  whether it is possible to rebuild the original CRT from its cut tree, how to proceed, and what additional information one needs to possess in order to make this happen. This question has been considered for the Brownian CRT \cite{BrWa} and for the stable L\'evy trees \cite{ADG}. However, the approaches in \cite{BrWa} and \cite{ADG} both rely heavily on the self-similar property of the underlying trees. These trees also happen to be compact.

In the present work, we address the question of reconstruction for models of continuum random trees that are not necessarily self-similar, nor compact. More precisely, let $(\cT, d, \mu)$ be a CRT belonging to the family of {\it Inhomogeneous continuum random trees} (ICRT) introduced by Aldous and Pitman \cite{Al00},  which are parametrised by $(\beta, \bth)$ where  
\begin{equation}\label{def: Theta}
	\beta\in [0, \infty), \quad \bth=(\theta_{1}, \theta_{2}, \dots) \quad\text{with } \  \theta_{1}\ge \theta_{2}\ge \cdots \ge 0  \ \text{ and } \ \sum_{i\ge 1}\theta_{i}^{2}<\infty. 
\end{equation}
We will assume that 
\begin{equation}
	\label{Theta-ass}
	\text{ if } \beta=0, \text{ then necessarily }\sum_{i\ge 1}\theta_{i}=\infty\,,
\end{equation}
to avoid trivial cases. 
Given such a pair $(\beta, \bth)$, the associated ICRT $(\cT, d, \mu)$ is constructed using the so-called Line-breaking Algorithm due to Aldous \& Pitman \cite{Al00}; we refer to Section \ref{sec:line_breaking} for details.   
Inhomogeneous continuum random trees are the scaling limits of birthday trees (or $\bp$-trees), a family of random trees which generalize uniformly random labelled trees and 
appear naturally in the study of the birthday problems \cite{Pi00}. 
They are also closely related to additive coalescents \cite{Al00} and models of inhomogeneous random graphs \cite{BHS}.

The Brownian CRT appears as a particular instance of the ICRT family by taking $\beta=1$ and $\theta_{i}=0$ for all $i\ge 1$. Whereas the Brownian CRT emerges in the limit of large uniform binary trees, as soon as $\bth$ is non zero, points of infinite degrees or {\it hubs} will be found in the corresponding ICRT. We denote by $\Br_{\infty}(\cT)$ the set of these hubs in $\cT$, which is necessarily a finite or countably infinite set. 
In the current work, we will focus on the case where $\sum_{i}\theta_{i}=\infty$. In that case, $\Br_{\infty}(\cT)$ is almost everywhere dense in $\cT$, and our approach to the reconstruction problem relies on a new approximation of the distance $d$ in $\cT$ using elements of $\Br_{\infty}(\cT)$. 
Let us also point out this distance approximation is used in \cite{Wa22+} to confirm a conjectured mixture relationship between stable L\'evy trees and ICRT. 

\medskip
\noindent
\textbf{Plan of the paper.} We will need a fair amount of preliminaries in order to be able to give the precise statement of our results. 
Section~\ref{sec:preliminaries} is therefore dedicated to the relevant background and results: In Section~\ref{sec:line_breaking}, we introduce the inhomogeneous continuum random trees via the Line-breaking Algorithm. We explain in Section~\ref{sec: approx} the new approach to distance approximations that is key to our reconstruction procedure. In Section~\ref{sec:pruning_cut-trees}, we define the pruning process on the ICRT and the accompanying cut tree. Our approach to the reconstruction problem is then described in Section~\ref{sec:reconstruction}. The rest of the paper contains the proofs of our main results. 

\section{ICRTs, fragmentation and cut trees}
\label{sec:preliminaries}

\subsection{The inhomogeneous continuum random trees} 
\label{sec:line_breaking}

Let $(\beta, \bth)$ be as in \eqref{def: Theta}. We will use the notation $\|\bth\|^{2}=\sum_{i\ge 1}\theta_{i}^{2}$. 
The following {\it Line-breaking} construction is a trivial extension to the original version presented in \cite{Al00}, where it is assumed that $\beta+\|\bth\|^{2}=1$. We build a (random) metric space $(\cT, d)$ with a collection of Poisson processes defined as follows. 
If $\beta>0$, let $\{(U_{j},W_{j}), j\ge 1\}$ be the atoms of a Poisson point measure on the first octant $\{(x,y) : 0\le y \le x\}$ with intensity $\beta$ per unit area. For every $\theta_{i}>0$, let $(\xi_{i, j} , j \ge 1)$ be the jumps of a Poisson process on the positive real line with intensity $\theta_{i}$ per unit length. Note that the hypotheses on $\bth$ entail that the only limit point of the set $\{U_{j}, j \ge 1, \xi_{i, j}, i \ge 1, j \ge 2\}$ is infinity, and therefore can be ordered as $0 < \eta_{1} < \eta_{2} < \dots$; we will refer to the $\eta_{i}$'s as \emph{cutpoints}. Next, we associate to each $\eta_{k}$ a \emph{jointpoint} $\eta_{k}^{\ast}$ in the following way. If $\eta_{k}=U_{j}$ for some $j\ge 1$, then $\eta_{k}^{\ast}=W_{j}$. If, instead, $\eta_{k}=\xi_{i, j}$ for some $j \ge 2$, we then set $\eta_{k}^{\ast} =\xi_{i, 1}$. Note that $\eta_{k}^{\ast}<\eta_{k}$ by definition for each $k$. 

Having got the cutpoints $(\eta_{k})_{k\ge 1}$, we use them to partition the half-line $[0, \infty)$ into line segments $[\eta_{k}, \eta_{k+1}]$, $k\ge 0$, with the understanding that $\eta_{0}=0$. We then assemble these line segments into a tree. More precisely, let $\cR_{1}$ be simply the single branch $[0, \eta_{1}]$. For $k\ge 1$, let $\cR_{k+1}$ be the (real) tree obtained from $\cR_{k}$ and $[\eta_{k}, \eta_{k+1}]$ by identifying $\eta_{k}$ with $\eta_{k}^{\ast}\in \cR_{k}$. Then $(\cR_{k})_{k\ge 1}$ is an increasing sequence of metric spaces. 
We then define $(\cT, d)$ to be the completion of $\cup_{k}\mathcal R_{k}$. 
It is often convenient to distinguish the image of $0$ in $\cT$; we refer to it as the {\it root} of $\cT$ and denote it as $\rho$. 
Let $\bP^{\beta, \bth}$ be the law of $(\cT, d, \rho)$. Note that it is straightforward to check from the properties of Poisson point process that for $c>0$, 
\begin{equation}
	\label{id: icrt-scale}
	(\cT, d) \text{ under } \bP^{c^{2}\beta,\, c\bth} \overset{(d)}{=} (\cT, \tfrac1c\, d) \text{ under }\bP^{\beta, \bth}.
\end{equation}

Say a point $\sigma\in \cT$ (resp.~$\sigma\in \cR_{k}$) is a {\it leaf} if $\cT\setminus\{\sigma\}$ (resp.~$\cR_{k}\setminus\{\sigma\}$) is still connected. It is not difficult to see that almost surely $\cR_{k}$ has exactly $k+1$ leaves, which are the images of $\eta_{0}, \eta_{1}, \eta_{2}, \dots, \eta_{k}$ in $\cR_{k}$. 
Let $\mu_{k}$ be the empirical measure of $\cR_{k}$ on the set of its leaves. It turns out  (\cite{Al00}) that $(\mu_{k})_{k\ge 1}$ converges in the weak topology of $\cT$ to a limit $\mu$, which is referred to as the 
{\it mass measure} of $\cT$. Moreover, $\mu$ has the following properties: it is non-atomic and supported on the set of leaves of $\cT$. 
Moreover, let $(V_{j})_{j\ge 0}$ be a sequence of i.i.d.~points sampled according to $\mu$ and let $\mathcal R'_{k}$ be the smallest subtree of $\cT$ containing $V_{0}, V_{1}, V_{2}, \cdots, V_{k}$ rooted at $V_{0}$; then $\mathcal R'_{k}$ has the same distribution as $\mathcal R_{k}$, the real tree obtained from the first $k$ branches in the previous construction. 
Note that $\cR_{k}$ can be seen as a deterministic function of $\cR_{k+1}$, and the same is true for $\cR'_{k}$ and $\cR'_{k+1}$. 
Hence, $(\cR'_{m}, \cR'_{k})\eqd (\cR_{m}, \cR_{k})$ for all $m\ge k\ge 1$. Letting $m\to\infty$, we deduce that for all $k\ge 1$, 
\begin{equation}
	\label{id: re-root}
	(\cT, \cR'_{k}) \text{ has the same law as } (\cT, \cR_{k}) \text{ under }\bP^{\beta, \bth}.
\end{equation} 

It also follows from the previous construction that for each $\theta_{i}>0$, the corresponding jointpoint $\xi_{i, 1}$ will be identified with the endpoints of an infinite number of branches. Denote by $B_i$ the image of $\xi_{i, 1}$ in $\cT$. 
For a real tree $(T, d)$, we define the {\it degree} of a point $\sigma\in T$ as the possibly infinite number of connected components of $T\setminus \{\sigma\}$, and write it as $\deg(\sigma, T)$. 
It will also be convenient to define $\deg(\sigma, T)=-\infty$ if $\sigma\notin T$. 
It is then not difficult to see that $\deg(B_{i}, \cR_{k})$ increases to infinity as $k$ grows. 
To measure its growth rate, we introduce
\begin{equation}
	\label{def: h}
	\Psi(t)=\tfrac12\beta t^{2}+\sum_{i\ge 1}(e^{-\theta_{i}t}-1+\theta_{i}t)\le \tfrac{t^{2}}{2}\Big(\beta + \|\bth\|^{2}\Big) <\infty, \quad t\ge 0.
\end{equation}
The assumptions \eqref{def: Theta} and \eqref{Theta-ass} ensure that $t\mapsto \Psi(t)$ is continuous, strictly increasing, and $\Psi(+\infty)=+\infty$. Let $\Psi^{-1}$ be its inverse function. 
Recall the set $\Br_{\infty}(\cT)$ of points of infinite degree in $\cT$. 
In Section \ref{sec: br-pt}, we prove:
\begin{prop}\label{prop:degrees}
	We have $\Br_{\infty}(\cT)=\{B_{i}: \theta_{i}>0, i\ge 1\}$ a.s. 
	For all $\sigma\in \cT$, the following limit exists in probability under $\bP^{\beta, \bth}$:
	\begin{equation}
		\label{def: br-pt}
		\Delta_{\cT}(\sigma):=\lim_{k\to\infty}\frac{\deg(\sigma, \cR_{k})}{\Psi^{-1}(k)}.
	\end{equation}
	Moreover, $\Delta_{\cT}(\sigma)>0$ if and only if $\sigma\in \Br_{\infty}(\cT)$; in that case we have $\Delta_{\cT}(B_{i})=\theta_{i}$, for each $B_{i}\in \Br_{\infty}(\cT)$. Let $(V_{i})_{i\ge 0}$ be a sequence of i.i.d.~points of $\cT$ with common law $\mu$, and let $\cR'_{k}=\cup_{1\le i\le k}\llb V_{0}, V_{k}\rrb$ be the subtree spanning $V_{0}, V_{1}, \dots, V_{k}$. Then 
	\[
	\big\{\deg(B_{i}, \cR'_{k}): k\ge 1, i\ge 1\big\} \text{ has the same law as } \big\{\deg(B_{i}, \cR_{k}): k\ge 1, i\ge 1\big\}\text{ under }\bP^{\beta, \bth}.
	\]
	As a consequence, we can replace $\cR_{k}$ by $\cR'_{k}$ in \eqref{def: br-pt} and obtain the same limit almost surely.
\end{prop}

We will refer to $\Delta_{\cT}(\sigma)$ as the {\it local time} of  $\sigma$ in $\cT$. 
Proposition~\ref{prop:degrees} implies that $\Delta_{\cT}(\sigma)$ is a measurable function in the Gromov--Prokhorov topology. 
It follows that almost surely
\[
\{\theta_{i}>0: i\in \N\}=\{\Delta_{\cT}(\sigma)>0: \sigma\in \cT\}=\{\Delta_{\cT}(\sigma): \sigma\in \Br_{\infty}(\cT)\} 
\]
is a measurable function of $(\cT, d, \mu)$. 

Let $\ell\in \N$ and suppose $\theta_{\ell}>0$. Let $\mathrm m^{(\ell)}=\sum_{1\le i\le \ell}\theta_{i}\delta_{B_{i}}$ be the finite measure on $\cT$ which puts a mass of size $\theta_{i}$ to the point $B_{i}$. We have for all $k\ge 1$, 
\begin{equation}
	\label{id: re-root'}
	\Big(\cR'_{k},\, \mathrm m^{(\ell)}(\cdot\cap \cR'_{k})\Big) \text{ has the same law as }\Big(\cR_{k}, \,\mathrm m^{(\ell)}(\cdot\cap \cR_{k})\Big) \text{ under }\bP^{\beta, \bth}.
\end{equation}
To see why this is true, we note that \eqref{id: re-root'} above is merely a reformulation of Proposition 5(b) in \cite{Al00}, in which $\cR_{k}$ was regarded as a discrete tree equipped with edge lengths and a subset of labeled vertices. 

\subsection{Approximations of the tree distance}
\label{sec: approx}

The following result will allow us to recover the length of a uniform branch in $\cT$ by tracking the local times 
of the points on that path. It is crucial to our approach for the reconstruction, which will essentially rely on identifying the points of infinite degree that used to be on any given path. 

\begin{prop}
	\label{prop: d-app} 
	Given the continuum random tree $(\cT, d, \mu)$, defined under $\bP^{\beta, \bth}$.
	Let $V$ and $V'$ be two independent points in $\cT$ with distribution $\mu$, let $\gamma_{\bth}(\epsilon)=\sum_{i}\theta_{i}\mathbf 1_{\{\theta_{i}>\epsilon\}}$, $\epsilon>0$. Suppose that $\sum_{i\ge 1}\theta_i=\infty$. Then under $\bP^{\beta, \bth}$, we have
	\begin{equation}\label{eq: icrt-app}
		\frac{1}{\gamma_{\bth}(\epsilon)}\sum_{b \in \llb V, V'\rrb \cap \Br_{\infty}(\cT)} \mathbf 1_{\{\Delta_{\cT}(b) > \epsilon \}} \xrightarrow{\epsilon\to 0} 
		d(V, V') \quad \text{ in probability.}
	\end{equation}
	If the entries of $\bth$ are all distinct, then the convergence in~\eqref{eq: icrt-app} also takes place almost surely. 
\end{prop}

As a consequence of Proposition~\ref{prop: d-app}, we can always find a sequence $\epsilon_{k}\to 0$, which only depends on $(\beta, \bth)$, so that for $\bP^{\beta,\bth}$-almost every realisation of $\cT$, we have
\begin{equation}
	\label{id: as-cv}
	d(V, V')=\lim_{k\to\infty}\frac{1}{\gamma_{\bth}(\epsilon_{k})}\sum_{b \in \llb V, V'\rrb \cap \Br_{\infty}(\cT)} \mathbf 1_{\{\Delta_{\cT}(b) > \epsilon_{k} \}} \quad \text{ in the a.s.~sense}.
\end{equation}
Proof of the above proposition will be given in Section~\ref{sec: icrt-app}. For the moment, let us simply note that the assumption $\sum_{i}\theta_{i}=\infty$ is always satisfied if $\beta=0$.

\subsection{Pruning continuum random trees} 
\label{sec:pruning_cut-trees}

Given the continuum random tree $(\cT, d, \mu)$, defined under $\bP^{\beta, \bth}$, 
we consider a pruning of $\cT$ according to a Poisson point measure. To that end, we first introduce $\mathcal L$, a $\sigma$-finite measure of $\cT$, which will serve as the intensity measure of this Poisson point measure. The {\it skeleton} of $(\cT,d)$ consists of the points of degree at least two and is defined as 
\[\Skel(\cT)= \bigcup_{\sigma,\sigma'\in \cT} \rrb \sigma, \sigma'\llb,\] 
where $\rrb \sigma,\sigma'\llb = \llb \sigma,\sigma'\rrb \setminus \{\sigma,\sigma'\}$ is the open arc between $\sigma$ and $\sigma'$.  
Let $\ell$ stand for the {\it length measure} of $\cT$, which is a $\sigma$-finite measure concentrated on $\Skel(\cT)$  
characterised by $\ell(\llb \sigma, \sigma'\rrb)=d(\sigma, \sigma')$ for all $\sigma, \sigma'\in \cT$. 
Recall also that to every point $b$ of infinite degree is associated a positive number $\Delta_{\cT}(b)$, as specified in Proposition~\ref{prop:degrees}. 
Recall the parameter $\beta$ from \eqref{def: Theta}. 
Let 
\begin{equation}\label{def: L}
	\mathcal L(dx) := \beta \cdot\ell(dx) +\sum_{b \in \Br_{\infty}(\cT)} \Delta_{\cT}(b)\cdot\delta_b (dx).
\end{equation}
We have $\mathcal L(\llb \sigma, \sigma'\rrb)\in (0, \infty)$, $\mathbf P^{\beta, \bth}$-a.s.~for all distinct $\sigma, \sigma'\in \cT$ (see Lemma~5.1 of \cite{BW}).

Now let $\mathcal P=\{(t_i, x_i): i\ge 1\}$ be a Poisson point measure on $\R_{+}\times\cT$ of  intensity  $dt\otimes\mathcal L(dx)$. An atom $(t_{i}, x_{i})$ of $\mathcal P$ is interpreted as a {\it cut} which disconnects the tree $\cT$ at location $x_{i}$ and at time $t_{i}$. So the pruning by $\mathcal P$ amounts to cutting the skeleton of the tree at a rate proportional to $\beta$ and removing each branch point $b$ of infinite degree at an exponential time of rate $\Delta_{\cT}(b)$. 
The above pruning process  on $\cT$ induces  a fragmentation process as follows. 
For $t\ge0$, denote by $\mathcal P_{t}:=\{x_i: \exists\, t_i \le t \text{ s.t. } (t_i, x_i)\in \mathcal P\}$ the set of the locations of cuts arriving before $t$; let us write $\cT\setminus\cP_{t}$ for the ``forest'' obtained from $\cT$ by removing the points in $\cP_{t}$. Let $\boldsymbol{\mu}_{t}=(\mu_{t}(i))_{i\ge 1}$ be the sequence of non zero $\mu$-masses of the components of $\cT\setminus\cP_{t}$ sorted in decreasing order. Alternatively, $\boldsymbol\mu_{t}$ can be recovered as the frequencies of an exchangeable random partition (\cite{Bert_frag}). More precisely, 
let $(V_i)_{i\ge 1}$ be i.i.d.~points of $\cT$ with common distribution $\mu$. 
For $t\ge 0$, say $i\sim_{t}j$ if and only if $i=j$ or $\llb V_{i}, V_{j}\rrb \cap \mathcal P_{t}=\varnothing$. 
Then the equivalence relation $\sim_{t}$ defines a partition of $\N$, denoted as $\Pi_{t}$. Moreover, we have a.s. 
\begin{equation}\label{def: Tt}
	i\sim_{t}j \quad \Longleftrightarrow\quad V_{j}\in \cT_{t}(V_{i}) :=\{\sigma\in \cT: \mathcal P_{t}\cap \llb V_{i}, \sigma\llb \,=\varnothing \}\,.
\end{equation}
For the block of $\Pi_{t}$ containing $i$, i.e.~$\{j\in \N: i\sim_t j\}$, its frequency, defined by
\[
\lim_{k\to\infty}\frac{1}{k}\Card\big\{j: 1\le j\le k, j\sim_{t} i\big\},
\]
exists a.s. by the Law of Large Numbers, and can be identified as $\mu(\cT_{t}(V_{i}))$. 
Moreover, $\boldsymbol\mu_{t}$ coincides with the rearrangement in decreasing order of the block frequencies of $\Pi_{t}$.

\subsection{Genealogy of the fragmentation and cut trees}
\label{sec: cut-tree}


The partition-valued process $(\Pi_{t})_{t\ge 0}$ has a natural genealogical structure: if $s>t$, then each block $B$ of $\Pi_{s}$ is contained in some block $B'$ of $\Pi_{t}$; in that case, we say that $B$ is a descendant of $B'$.   
It turns out that this genealogical structure can be encoded by another CRT called the \emph{cut tree}, which also contains information on the times of fragmentation events. 
This has been explored in \cite{BeMi13, BW, Di15, AbDe13}. Let us briefly recall the construction of the cut tree.

Write $\mathbb N_{0}:=\{\mathrm{0}, \mathrm{1}, \mathrm{2}, \dots\}$. For $i\ne j$, let 
\begin{equation}
	\tau_{ij}=\inf\{t>0: \cP_{t}\cap\llb V_{i}, V_{j}\rrb \ne \varnothing\}
\end{equation}
be the first time when a cut falls on $\llb V_i, V_j\rrb$, which is distributed as an exponential variable of rate $\mathcal L(\llb V_{i}, V_{j}\rrb)$. We introduce a function $\delta: \N_{0}\times \N_{0} \to [0, \infty]$ whose values are given as follows. Set $\delta(i, i)=0$ for all $i\in \N_{0}$; let
\begin{equation}
	\label{def: delta}
	\delta(0, i) = \delta(i, 0)= \int_0^\infty \mu\big(\cT_{s}(V_{i})\big)ds\in [0, +\infty], \quad i\in \N,
\end{equation}
where $\cT_{s}(V_{i})$ is the part of $\cT$ still connected to $V_{i}$ at time $s$, as defined in \eqref{def: Tt};  for $i, j\ge 1$ and  $i\ne j$, define
\begin{equation}
	\label{def: delta'}
	\delta(i, j) = \delta(j, i) = \int_{\tau_{ij}}^{+\infty} \Big\{\mu\big(\cT_{s}(V_{i})\big)+\mu\big(\cT_{s}(V_{j})\big)\Big\} ds. 
\end{equation}

\begin{thm}[\cite{BW}, Theorem 3.5]
	\label{prop: c_dist}
	$(\delta(i, j))_{i, j\in \N_0}$ has the same distribution as  $(d(V_{i+1}, V_{j+1}))_{i, j\in \N_{0}}$  under $\bP^{\beta, \bth}$. 
\end{thm}

In particular, Theorem~\ref{prop: c_dist} implies that $\delta(i,j)<\infty$ almost surely for all $i, j\in \N_{0}$. Furthermore, it can be readily checked from the definitions \eqref{def: delta} and \eqref{def: delta'} that $\delta$ is a metric on $\N_{0}$ which verifies the four-point condition, namely $\delta(i, j)+\delta(k, l)\le \max\{\delta(i, k)+\delta(j, l), \delta(j, k)+\delta(i, l)\}$, for any $i,j,k,l\in \N_0$. 
Let $(\cC, \delta)$ be the metric completion of $(\N_{0}, \delta)$. It turns out that $(\cC, \delta)$ is connected, and thus a real tree. We define $0$ to be its root. Denote by $\mathcal S_{k}=\cup_{1\le i\le k}\llb 0, i\rrb$ the subtree of $\cC$ spanning $\{0, 1, 2, \dots, k\}$ and by $\mathcal R'_{k}=\cup_{0\le i\le k}\llb V_{0}, V_{i}\rrb$ the subtree of $\cT$ spanning $\{V_{0}, V_1,\dots, V_{k}\}$, which we see as real trees equipped with their respective metrics. It follows from Theorem~\ref{prop: c_dist}  that $\mathcal S_{k}$ has the same distribution as $\mathcal R'_{k}$, for each $k\ge 1$. In particular, $(\mathcal S_{k})_{k\ge 1}$ satisfies the so-called {\em leaf-tight} property: 
\[\inf_{i\ge 2} \delta(1,i) = 0 \quad \text{a.s.}\]
See Lemma 6 in \cite{Al00} for a proof of this. Appealing to the general theory in \cite{Aldous1993a}, and the arguments in \cite{BeMi13}, we see that the empirical measures $\frac{1}{k}\sum_{1\le i\le k}\delta_{i}$, converge weakly to some limit $\nu$. 
In fact, we can rephrase the distributional identity in Theorem \ref{prop: c_dist} as follows. 

\smallskip
\noindent
\textbf{Theorem \ref{prop: c_dist}* }(\cite{BW}){\bf. }
The measured metric space $(\cC, \delta, \nu)$ is distributed as $(\cT, d, \mu)$ under $\bP^{\beta, \theta}$. 

\smallskip

Let us point out two immediate consequences of the above construction of $(\cC, \delta, \nu)$: (i)  $(\cC, \delta, \nu)$ is a measurable function of $(\cT, d, \mu)$ and $\cP$, and (ii) $\N$ is an i.i.d.~sequence of points with distribution $\nu$, in particular the whereabouts of the points in $\N$ are not retained in $\cC$. As we will see, perhaps unsurprisingly, in order to recover the genealogy of the fragmentation, one  needs the locations of $\N$ in the cut tree. 

\begin{figure}[tp]
	\centering
	\includegraphics[height = 5cm]{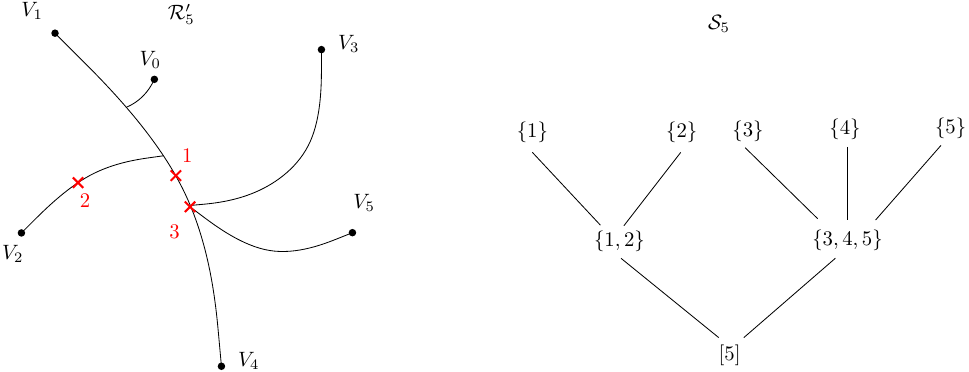}
	\caption{\label{fig1}An example of the pruning restricted to $\mathcal R'_{k}$ and the accompanying $\mathcal S_{k}$. On the left, an illustration of the spanning tree $\cR'_{k}$ with $k=5$. Locations of cuts are depicted as $\times$, with the numbers next to them indicating their arrival orders. On the right, the shape of $\mathcal S_{k}$ which captures how $[k]$ splits into singletons due to the cuts on the left-hand side.}
\end{figure}

Now given the pair $(\cC, \N)$, we explain how to recover the genealogy of $(\Pi_{t})_{t\ge 0}$. To that end, let us first take $k\ge 2$ and consider the restriction of the pruning to the spanning tree $\mathcal R'_{k}=\cup_{1\le i\le k}\llb V_{0}, V_{i}\rrb$. Depending on whether $V_{i}$ is connected to $V_{j}$, $1\le i, j\le k$, at time $t$ in this pruning on $\cR'_{k}$, we obtain a partition $\Pi^{[k]}_{t}$ of the set $[k]:=\{1, 2, \dots, k\}$. Note that $\Pi^{[k]}_{t}$ is simply the restriction of $\Pi_{t}$ to $[k]$. Moreover, the aforementioned ancestor-descendant relation reduces to as follows. The root of  $\cS_{k}=\cup_{1\le i\le k}\llb 0, i\rrb$ corresponds to the common ancestor $[k]=\Pi^{[k]}_{0}$. The first time a cut falls on $\mathcal R'_{k}$, $[k]$ splits into two or more blocks, and these blocks form the first generation in $\cS_{k}$. 
As more cuts arrive on $\mathcal R'_{k}$, each of these blocks splits until $\Pi^{[k]}_{t}$ eventually is comprised solely 
of singletons, which form the last generation. In other words, the genealogy induced by this ancestor-descendant relation is captured by the 
shape of $\mathcal S_{k}$ (see Fig.~\ref{fig1}).  
Letting $\cS_{k}$ grow to $\cC$, we see that the entire genealogy of $(\Pi_{t})_{t\ge 0}$ is retained in $(\cC, \N)$. 

\section{Main results on the reconstruction from cut trees} 
\label{sec:reconstruction}

Our aim here is to devise a way to rebuild the initial tree $(\cT, d, \mu)$ from a possibly augmented version of the cut tree $(\cC, \delta, \nu)$. 
We begin with an informal description of our strategy. 
Let us recall the i.i.d.~sequence $(V_{i})_{i\ge 1}$ of random points with distribution $\mu$ on $\cT$. The CRT theory (\cite{Al00}) says that $(\cT, d, \mu)$ is itself a measurable function of the semi-infinite matrix $(d(V_{i}, V_{j}))_{i, j\ge 1}$. Thus, as soon as the latter is recovered, the continuum random tree will follow suit. 

As Proposition~\ref{prop: d-app} suggests, one way to recover $d(V_{i}, V_{j})$ is to find out which points of infinite degree lie on the path $\llb V_{i}, V_{j}\rrb$ and what their local times are. 
To answer these two questions, we get a helpful hand from Proposition~\ref{prop: br_pt} (see below), which says that there is a one-to-one correspondence between $\Br_{\infty}(\cC)$ and $\Br_{\infty}(\cT)$ that preserves the local times, where $\Br_{\infty}(\cC)$ and $\Br_{\infty}(\cT)$ are the respective sets of points of infinite degree in $\cC$ and $\cT$. Call $\phi$ this map from $\Br_{\infty}(\cC)$ to $\Br_{\infty}(\cT)$; we will be able to read the value of $\Delta_{\cT}(\phi(b))$ as $\Delta_{\cC}(b)$ for every $b\in \Br_{\infty}(\cC)$.

It remains to see how we can locate the set $\phi^{-1}(\Br_{\infty}(\cT)\cap \llb V_{1}, V_{2}\rrb)$ in $\Br_{\infty}(\cC)$. As it turns out, this piece of information is only partially encoded in $\cC$. Therefore, to enable an a.s.~reconstruction of the initial tree, we decorate the cut tree $\cC$ with some additional random variables called {\it traces}.  To explain this idea of traces, let us first take a look at the reconstruction problem for discrete trees, or more precisely, a model of random trees called $\bp$-trees, which can be seen as the discrete counterparts of ICRT.

\subsection{Reconstruction of $\bp$-trees}
\label{sec: intro-p-tree}

For a probability measure $\bp$ of the set $[n]=\{1, 2, 3, \dots, n\}$, let us write $p_{i}$ for its mass at $i$. It will be convenient to assume that $p_{1}\ge p_{2}\ge \cdots p_{n}>0$. 
We regard a rooted tree as a family tree, with the root as the common ancestor and the neighbours of the root as the first generation, etc. 
Denote by $\mathbb T_{n}$ the set of all rooted trees with vertices labelled by $[n]$. For each $t\in \mathbb T_{n}$, we define
\begin{equation}
	\label{def: pi}
	\pi_{\bp}(t)=\prod_{i\in [n]}p_{i}^{D^{+}_{i}(t)}, 
\end{equation}
where $D^{+}_{i}(t)$ is the number of children of vertex $i$ in $t$. Cayley's multinomial formula \cite{Pi99} says that 
\begin{equation}
	\label{eq: Cayley}
	\sum_{t\in \mathbb T_{n}}\pi_{\bp}(t)=\Big(\sum_{i} p_{i}\Big)^{n-1}=1\,,
\end{equation}
and therefore, $\pi_{\bp}$ defines a probability measure on $\mathbb T_{n}$. 
A random element of $\mathbb T_{n}$ with distribution $\pi_{\bp}$ will be referred to as a $\bp$-tree. 

\medskip
\noindent\textbf{Pruning $\bp$-trees and discrete cut trees.} 
Let $T$ be a $\bp$-tree. We introduce a pruning procedure of $T$ by removing vertex $i$ 
together with all the edges adjacent to it 
at a random time $E_{i}$. We further assume that $E_{i}$, $1\le i\le n$, are independent exponential variables with 
$\mathbb E[E_{i}]=1/p_{i}$.

\begin{figure}[tp]
	\centering
	\includegraphics[height = 5.5cm]{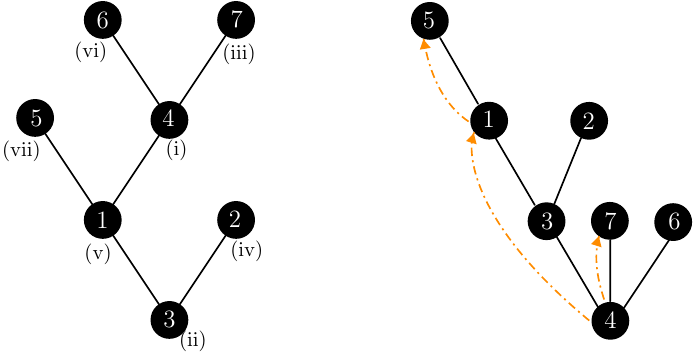}
	\caption{\label{fig2} An example of pruning of $\bp$-tree and its cut tree. On the left, a realisation of the tree $T$ with vertices labeled from 1 to 7. The roman numbers next to the vertices refer to the order in which they are removed in the pruning. On the right, the corresponding cut tree $C$. Observe that Vertex 1 has degree 3 on the left and degree 2 on the right. In this example, Vertex 4 has three subtrees above it in $C$; the traces associated to the three subtrees, from left to right,  are respectively 1, 7, 6. To reconstruct the path of $T$ from 5 to 7, we note that their most recent common ancestor in $C$ is 4.  Since 7 is a trace of 4, $\{4, 7\}$ is an edge of $T$. Similarly, $\{4, 1\}$ is also an edge of $T$. Since 5 is a trace of 1, we can move from 5 to 7 in $T$ by following the edges $\{5, 1\}, \{1, 4\}$ and $\{4, 7\}$.  }
\end{figure}

We define here another rooted tree $C$, which will serve as the discrete analogue of the cut tree and encodes 
information including the genealogy of the above pruning process. 
Note that our definition is different from the one in \cite{BeMi13}, but ensures our cut tree $C$ is defined on the same vertex set as $T$.  
The first removed vertex, say $X_{1}$, is the root of $C$. 
Removing $X_{1}$ splits the tree $T$ into a number of connected components, say $R_{1}, R_{2}, \dots, R_{m}$. The vertex sets of 
these components induce a partition of $[n]\setminus\{X_{1}\}$. 
Similarly, removing the root $X_{1}$ from $C$ results in a collection of subtrees, say $S_{1}, S_{2}, \dots, S_{m'}$. The tree $C$ is set up in such a way that the vertex sets of these subtrees coincide with 
those of $R_{1}, R_{2}, \dots, R_{m}$. 
In particular, this implies $m'=m$. 
We then iterate this procedure for each of the remaining components $(R_{i})_{1\le i\le m}$, so that $S_{i}$ is the cut tree of $R_{i}$. 
More precisely, the root of $S_{i}$ is the first removed vertex in $R_{i}$, and so on so forth. 
See  Fig.~\ref{fig2} for an example. 

\medskip
\noindent\textbf{Rebuilding $T$ from $C$.} 
Rebuilding $T$ amounts to identifying its edge set. To that end, we observe that if $v\in [n]$ has $m$ children in $C$, then $m$ edges in total are removed exactly at the time when $v$ is removed in the pruning. 
Denoting $D_{v}(T)$ for the degree of $v$ in $T$ (so that $D_{v}(T)-1=D^{+}_{v}(T)$ apart from the root), we have
\begin{equation}
	\label{ineq-degree}
	D_{v}(T)\ge m = D^{+}_{v}(C), 
\end{equation}
and we note that the inequality can be strict because of the edges removed in the previous steps. 
To recover these $m$ edges, we make the following observation: say $\mathcal V_{i}$, $1\le i\le m$, are the vertex sets of the subtrees of $C$ rooted at the children of $v$; then for $1\le i\le m$, there is a unique edge $e$ in $T$ between $v$ and a vertex $u_{i}\in \mathcal V_{i}$. 
We call $(u_{i})_{1\le i\le m}$ the {\it traces} for $v$ in $C$, which are randomly distributed vertices conditioned on $C$. 
Observing that there is a unique trace for $v$ in each of the subtrees above $v$, we define the general collection of traces as follows.
For each vertex $v\in [n]$, let $C_{v, 1}, C_{v, 2}, \dots, C_{v, D^{+}_{v}(C)}$ be the subtrees of $C$ above $v$, ranked in the decreasing order of their $\bp$-masses (breaking ties arbitrarily).  
Let $u_{v, i}$ be the trace for $v$ that belongs to $C_{v, i}$, $1\le i\le D^{+}_{v}(C)$. 
We can reconstruct $T$ from $C$ as soon as we have the knowledge of the collection of traces $\{u_{v, i}: 1\le i\le D^{+}_{v}(C), v\in [n]\}$. Results from \cite{BW} say that this collection consists of independent points, with each $u_{v, i}$ distributed according to $\bp$ conditioned to be in $C_{v, i}$ (see Lemma~\ref{lem: p-trace} below).

There is another approach to this notion of trace that will be useful to us later on. Fix two distinct vertices $w_{1}, w_{2}\in [n]$. Knowing the collection of traces  $\{u_{v, i}: 1\le i\le D^{+}_{v}(C), v\in [n]\}$ will allow us to reconstruct the path in $T$ linking $w_{1}$ to $w_{2}$, denoted as $\llb w_{1}, w_{2}\rrb$.  Indeed, let $b_{0}=w_{1}\wedge w_{2}$ be the most recent common ancestor of $w_{1}$ and $w_{2}$ in $C$. Then $b_{0}$ is in fact the location of the first cut falling on the path $\llb w_{1}, w_{2}\rrb$ during the pruning of $T$. Assume that $b_{0}$ is distinct from $w_{1}$ and $w_{2}$; then there are two subtrees of $C$ above $b_{0}$ containing $w_{1}$ and $w_{2}$ respectively. Let $u_{1}$ and $u_{2}$ be the traces for $b_{0}$ contained in these two subtrees. Then it should be clear from the previous paragraph that $u_{1}, u_{2}$ are the two neighbours of $b_{0}$ in the path $\llb w_{1}, w_{2}\rrb$. If $u_{1}\ne w_{1}$ or $u_{2}\ne w_{2}$, then we can repeat this procedure on the subpath $\llb u_{1}, w_{1}\rrb$ or $\llb u_{2}, w_{2}\rrb$, and so on so forth. This leads to the identification of a sub-collection of the traces which we will refer to as the routings of the pair $(w_{1}, w_{2})$. Following the most recent common ancestors in $C$ of these routings will allow us to locate all the cuts on $\llb w_{1}, w_{2}\rrb$, which will provide enough information for recovering the path $\llb w_{1}, w_{2}\rrb$. See Fig.~\ref{fig2} for an example.

\subsection{Images in $\cC$ of cuts in $\cT$}
\label{sec: phi-def}

Let us recall the set $\cP_{t}=\{x_{i}: \exists\, t_{i}\le t \text{ such that } (t_{i}, x_{i})\le \cP\}$ and denote by $\cP_{\infty} =\cup_{t\ge 0}\cP_{t}$. 
Let $\Br(\cC)=\{\sigma\in \cC: \deg(\sigma, \cC)\ge 3\}$ be the set of branch points in $\cC$ and we define similarly the set $\Br(\cT)$. In particular,  $\Br_{\infty}(\cC)$, the set of infinite degree points of $\cC$, is a subset of $\Br(\cC)$, and similarly $\Br_{\infty}(\cT)\subseteq \Br(\cT)$. 
There is a natural correspondence between $\Br(\cC)$ and $\cP_{\infty}$ that we now explain. 
Recall from Section~\ref{sec: cut-tree} the way $\cC$ is defined and how it is related to the genealogy of the partition-valued process $(\Pi_{t})$. 
Let $i\wedge j$ be the most recent common ancestor of $i$ and $j$ in $\cC$; this is where the two paths from the root  leading respectively to $i$ and $j$ separate, and corresponds to the event at time $\tau_{i, j}$ when $V_i$ and $V_j$ cease to belong to the same connected component.
We define a map 
\[
\phi: \Br(\cC) \to \cP_{\infty}
\]
by setting $\phi(i\wedge j) = x_{i, j}$, where $x_{i, j}\in \cP_{\infty}$ is the unique point of $\cT$ satisfying $(\tau_{i, j}, x_{i, j})\in \mathcal P$. Note that the above definition does not depend on the particular choices of $i, j$.
Indeed, if $(m, n)\in \N^{2}$ also satisfies $m\wedge n=i\wedge j$, then it can be checked that  we have $\tau_{i, j}=\tau_{m, n}$. 
Since all elements of $\Br(\cC)$ can be written as $i\wedge j$ for some $i, j\in \N$, $\phi$ is indeed a function defined on $\Br(\cC)$. 

We recall from \eqref{def: L} that the intensity measure $\cL$ of $\cP$ does not charge the set of binary branch points of $\cT$ nor the leaf set. It follows that if $x\in \cP_{\infty}$, then $\deg(x, \cT)\in \{2, \infty\}$ almost surely. 
Since $\cC$ is an ICRT itself (Theorem~\ref{prop: c_dist}), the local time $\Delta_{\cC}(\sigma)$ exists for every branch point $\sigma\in \Br(\cC)$. It turns out that the map $\phi$ has the following properties.

\begin{prop}\label{prop: br_pt}
	Let $\phi: \Br(\cC)\to \cP_{\infty}$ be as defined above. The following statements hold 
	$\mathbf P^{\beta, \bth}$-a.s. 
	\begin{compactenum}[i)]
		\item If $\sigma\in \Br(\cC)\setminus\Br_{\infty}(\cC)$, then $\phi(\sigma)$ has degree 2 in $\cT$;
		\item If $\sigma \in \Br_\infty(\cC)$, then $\phi(\sigma)\in \Br_{\infty}(\cT)$ and $\Delta_{\cC}(\sigma)=\Delta_{\cT}(\phi(\sigma))$. 
	\end{compactenum}
\end{prop}

Proposition~\ref{prop: br_pt} is proved in Section~\ref{sec:proof_correspondence}. For now, let us simply point out that for the discrete tree $T$ and its cut tree $C$ introduced in Section~\ref{sec: intro-p-tree}, the role of $\phi$ is played by the identity map between the vertex sets of the two trees. However, when it comes to comparing the vertex degrees, we only have the inequality \eqref{ineq-degree} in general.

\subsection{Traces and routings in the ICRT}
\label{sec: trace-def}

In this part, we introduce the notion of trace for inhomogeneous continuum random trees. This is a countable collection of points in the cut tree $\cC$. Together with $\N$, the two collections represent the additional data required for the reconstruction of $\cT$. The idea of trace has already been introduced for compact continuum random trees in \cite{ADG}, albeit under the term {\it image}. We extend it here to complete continuum random trees which are not necessarily compact.  

To motivate the definition of trace and routing in the continuum random tree setting, let us recall from Section~\ref{sec: intro-p-tree} that in the discrete cut trees, the traces for the vertex $v\in [n]$ are a collection of vertices, each belonging to a different subtree above $v$. 
Now for each $b\in \Br(\cC)$, write $\cC_{b}=\{\sigma\in \cC: b\in \llb 0, \sigma\llb\}$ to denote the {\it subtree of $\cC$ above $b$}.  
Define $\mathrm D_{b}=\{1, 2\}$ if $b$ is a binary branch point, and $\mathrm D_{b}=\N$ if $b\in \Br_{\infty}(\cC)$.  
Denote by $(\cC_{b, i})_{i\in \mathrm D_{b}}$ the collection of connected components of $\cC_{b}$ ranked in decreasing order of their $\nu$-masses. A {\it collection of traces} in $\cC$ is a collection $\{\sigma_{b, i}: i\in \mathrm D_{b}, b\in \Br(\cC)\}$ that satisfies $\sigma_{b, i}\in \cC_{b, i}$ for all $i$ and $b$. 

\begin{figure}
	\centering
	\includegraphics[height = 7cm]{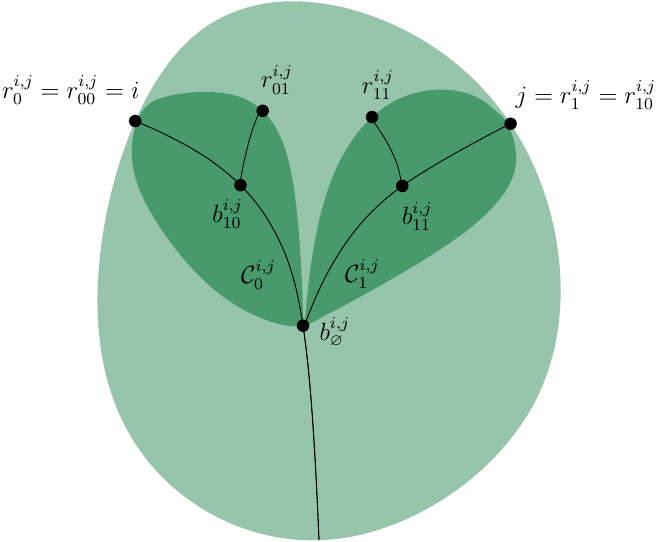}
	\caption{\label{fig3} An illustration of the first 2 generations in $\mathcal R_{W}^{i, j}$.}
\end{figure}

Suppose that $W$ is a collection of traces. 
Let $i, j\in\N$ be distinct. We define a sub-collection $\mathcal R_{W}(i, j)$ of $W$, called the {\it routing for the pair $(i, j)$ associated to $W$}, by taking inspiration from the previous reconstruction of a path in the discrete tree; see Fig.~\ref{fig3} for an illustration. 
Write $\mathbb U=\cup_{n\ge 0}\{0, 1\}^{n}$ for the complete binary tree and $\mathbb U^{\ast}=\mathbb U\setminus\{\varnothing\}$. 
Set $r^{i, j}_{0}= i$ and $r^{i, j}_{1}=j$. Suppose that $r^{i, j}_{u}$ has been defined for all $u\in \{0, 1\}^{n}$, $n\ge 1$. 
Then set $r^{i, j}_{u0}=r^{i, j}_{u}$. 
Let $w\in \{0, 1\}^{n-1}$ be the prefix of $u$, i.e.~$u=wq$ for some $q\in \{0, 1\}$. Assume $r^{i, j}_{w0}\ne r^{i, j}_{w1}$; then  their most recent common ancestor in $\cC$ is well-defined, denoted as $b^{i, j}_{w}=r^{i, j}_{w0}\wedge r^{i, j}_{w1}$. 
Write $\cC^{i, j}_{w0}$ (resp.~$\cC^{i, j}_{w1}$) for the connected component of $\cC_{b^{i, j}_{w}}$ that contains $r^{i, j}_{w0}$ (resp.~$r^{i, j}_{w1}$). 
Note that $\cC^{i, j}_{w0}$ must be somewhere in the list $(\cC_{b, i})_{i\in \mathrm D_{b}}$ for $b=b^{i, j}_{w}$. Namely, 
there is some $k\in \mathrm D_{b}$ so that $\cC^{i, j}_{w0}=\cC_{b, k}$. 
Let $r^{i, j}_{w01}=\sigma_{b, k}$ be the element of $W$ that is associated to the subtree $\cC^{i, j}_{w0}$. 
Similarly, let $k'\in \mathrm D_{b}$ be such that $\cC^{i, j}_{w1}=\cC_{b, k'}$. 
Let $r^{i, j}_{w11}=\sigma_{b, k'}$ be the element of $W$ that is associated to the subtree $\cC^{i, j}_{w1}$. 
We have thus defined  $r^{i, j}_{u}$ for all $u\in \{0, 1\}^{n+1}$. 
Inductively, this gives us the routing for $(i, j)$: $\mathcal R_{W}(i, j)=\{r^{i, j}_{u}: u\in \mathbb U^{\ast}\}$. 
We also denote $\mathcal B_{W}(i, j)=\{r^{i, j}_{u0}\wedge r^{i, j}_{u1}: u\in \mathbb U\}$. 

\begin{prop}
\label{prop: route}
There is a collection of traces $\mathcal W=\{\sigma_{b, i}: i\in \mathrm D_{b}, b\in \Br(\cC)\}$ so that the associated routings $\mathcal R_{\mathcal W}(i, j)$, $\mathcal B_{\mathcal W}(i, j)$ are well-defined for all $i\ne j$. Let $\phi: \Br(\cC)\to \cP_{\infty}$ be the map in Proposition~\ref{prop: br_pt}. We have
\[
\phi(b) \in \llb V_{i}, V_{j}\rrb \; \Longleftrightarrow \; b \in \mathcal B_{\mathcal W}(i, j). 
\]
\end{prop}

\begin{prop}[Distribution of the traces]
\label{prop: tr-dist}
Let $\mathcal W=\{\sigma_{b, i}: i\in \mathrm D_{b}, b\in \Br(\cC)\}$ be the collection of traces in Proposition~\ref{prop: route}. Given $\cC$, the collection $\mathcal W$ is independent, and $\sigma_{b, i}$ has the law $\nu$ restricted to $\cC_{b, i}$ for all $i\in \mathrm D_{b}$ and $b\in \Br(\cC)$. 
\end{prop}

The proofs of these two propositions are given in Section~\ref{sec: trace-pf}.

\subsection{Reconstruction of the ICRT}

Recall from \eqref{def: Theta} the parametrisation  of the ICRT family. Throughout this subsection, we assume that 
\begin{equation}
\label{hyp: theta}
\sum_{i\ge 1}\theta_{i}=\infty.
\end{equation}
The alternative case is somehow simpler, and will be discussed in Section~\ref{sec: discuss}. 

Let $\mathcal W$ be the collection of traces in Proposition~\ref{prop: route}, and let $\mathcal R_{\mathcal W}(i, j), \mathcal B_{\mathcal W}(i, j)$ be the associated routings, $i\ne j$. Combining Propositions~\ref{prop: d-app}, \ref{prop: br_pt} and \ref{prop: route} together immediately leads to the following

\begin{thm}[Almost sure reconstruction]
Let $\gamma_{\bth}(\epsilon)=\sum_{i}\theta_{i}\mathbf 1_{\{\theta_{i}>\epsilon\}}$ and assume \eqref{hyp: theta} holds. Let $(\epsilon_{k})_{k\ge 1}$ be the sequence in \eqref{id: as-cv}. We have $\bP^{\beta, \bth}$-a.s.
\begin{equation}
	\label{eq: recons}
	d(V_{i}, V_{j}) = \lim_{k\to\infty}\frac{1}{\gamma_{\bth}(\epsilon_{k})}\sum_{b\in \mathcal B_{\mathcal W}(i, j)}\mathbf 1_{\{\Delta_{\cC}(b)> \epsilon_{k}\}}  
\end{equation}
for all $i, j\in \N$. 
\end{thm}

Since $(\cT, d, \mu)$ can be recovered from $(d(V_{i}, V_{j}))_{i, j\ge 1}$, the above implies that $(\cT, d, \mu)$ is a measurable function of $(\cC, \N, \mathcal W)$ in the Gromov--Prokhorov topology (see for instance \cite{GPW09} for a definition). On the other hand, when only $(\cC, \delta, \nu)$ is given, Proposition~\ref{prop: tr-dist} indicates a way to re-sample the missing data.  Indeed, 
let $\tilde\N$ be a sequence of i.i.d.~points with common law $\nu$.
Let also $\tilde{\mathcal W}=\{\tilde\sigma_{b, i}: i\in \mathrm D_{b}, b\in \Br(\cC)\}$ be an independent collection where $\tilde\sigma_{b, i}$ is distributed according to $\nu$ conditioned to be in $\cC_{b, i}$. Then we have $(\cC, \tilde\N, \tilde{\mathcal W})\eqd (\cC, \N, \mathcal W)$. 
We then proceed to define the associated routings $\tilde{\mathcal R}_{\tilde{\mathcal W}}(i, j)=\{\tilde r^{i, j}_{u}: u\in \mathbb U^{\ast}\}$ by following the steps previously described. 
Since $\tilde\sigma_{b, i}$ are independently sampled, it can be checked with an inductive argument that $\tilde r^{i, j}_{w0}\ne \tilde r^{i, j}_{w1}$ for all $w\in \mathbb U$. Therefore, the routings are all well-defined, and we can carry out a reconstruction of $\cT$ by using $\tilde{\mathcal R}_{\tilde{\mathcal W}}(i, j), i, j\in \N$, as the input.
Denote by $\tilde\cT$ the obtained CRT. We then have $(\cC, \tilde\N, \tilde{\mathcal W}, \tilde\cT)\eqd ( \cC, \N, \mathcal W, \cT)$, and consequently $(\cC, \tilde\cT)\eqd (\cC, \cT)$.

\subsection{Discussion}
\label{sec: discuss}

\noindent
{\bf Reconstruction of ICRT with a Brownian component. } When \eqref{hyp: theta} fails, one must have $\beta>0$ by assumption. The reconstruction procedure in this case is actually much simpler, and we explain it right away: assume $\beta>0$ and take $\sigma, \sigma'\in \cT$. Let us denote by $\mathrm N_{t}(\sigma, \sigma')=\#(\cP_{t}\cap \llb \sigma, \sigma'\rrb\setminus \Br(\cT))$ the number of cuts that have fallen on $\llb \sigma, \sigma'\rrb \setminus \Br(\cT)$ before time $t$. We note that $\mathrm N_{t}(\sigma, \sigma'), t\ge 0$, is a Poisson process with rate $\beta d(\sigma, \sigma')$. It then follows from the Law of Large Numbers  that
\[
\frac{1}{\beta t} \mathrm N_{t}(\sigma, \sigma') \xrightarrow[]{t\to\infty} d(\sigma, \sigma') \quad \text{a.s.}
\]
Let us now explain how to determine $\mathrm N_{t}(\sigma, \sigma')$ from $(\cC, \N, \mathcal W)$. Recall that for each $b\in \Br(\cC)$, there is a unique time $\tau(b)$ so that $(\tau(b), \phi(b))\in \cP$. Indeed, if $b=i\wedge j$, then we have $\tau(b)=\tau_{i, j}$ the moment when $V_{i}$ is separated from $V_{j}$. As shown in Theorem 16 of \cite{ADG} (whose proof works for general CRT), $\tau(b)$ can be expressed as follows. Let $\ell_{\cC}$ stand for the length measure of $\cC$ and let $\cC_{x}=\{\sigma\in \cC: x\in \llb 0, \sigma\rrb\}$ denote the subtree of $\cC$ above $x\in\cC$. We have
\[
\tau(b)=\int_{\llb 0, b\rrb}\frac{d\ell_{\cC}(x)}{\nu(\cC_{x})}, \quad \text{a.s.}
\]
Denote by $\Br_{2}(\cC)=\{\sigma\in \cC: \deg(\sigma, \cC)=3\}$ the set of all binary branch points in $\cC$. We note that 
if $\phi(b)\in \cP_{t}\cap \llb \sigma, \sigma'\rrb\setminus \Br(\cT)$, then according to Proposition~\ref{prop: br_pt}, we have $b\in \Br_{2}(\cC)$. 
Combined with the previous arguments, this implies that for all $i, j\in \N$, 
\[
d(V_{i}, V_{j})=\lim_{t\to\infty}\frac{1}{\beta t} \sum_{b\in \mathcal B_{\mathcal W}(i, j)\cap \Br_{2}(\cC)} \mathbf 1_{\{\tau(b)\le t\}} \quad \text{a.s.}
\]

\medskip
\noindent
{\bf Reconstruction of L\'evy trees.} 
We believe that our current approach to the reconstruction problem can also apply to another family of continuum random trees: L\'evy trees, which will include the case of stable trees treated in \cite{ADG}. We discuss here what this will involve. Let $\alpha, \beta\ge 0$ and $\pi$ be a $\sigma$-finite measure on $(0, \infty)$ satisfying $\int_{(0, \infty)}(x\wedge x^{2})\pi(dx)<\infty$. We define a function $\Psi: \R_{+}\to \R_{+}$ as follows:  
\[
\Psi(\lambda)=\alpha\lambda+\tfrac{1}{2}\beta t^2+\int_{(0, \infty)}\big(e^{-\lambda x}-1+\lambda x\big)\pi(dx), \quad \lambda\ge 0,
\]
We further assume that $\int^\infty d\lambda/\Psi(\lambda)<\infty$, which holds true if $\beta>0$ or $\pi((x, \infty))$ is regularly varying at $0+$. 
For such a  function $\Psi$, it is then possible to define a continuum random tree $\cT$, called the $\Psi$-L\'evy tree, under the so-called excursion measure $\mathbf N^{\Psi}$; see \cite{DuLG02, DuLeG05}.  Roughly speaking, $\cT$ under $\mathbf N^{\Psi}$ represents the genealogy, issued from a single individual,  which evolves as a continuous-state branching process with branching mechanism $\Psi$. 

The L\'evy tree $\cT$ shares some strikingly similar features with inhomogeneous continuum random trees. In particular, as soon as $\pi$ is non zero, $\cT$ contains branch points of infinite degree. For each such branch point $b$, it is also possible to introduce an analogue notion of  local time  $\Delta_{\cT}(b)$, although the collection $\{\Delta_{\cT}(b): b\in \Br(\cT)\}$ is no longer deterministic in general and relates to jumps in a L\'evy process. As a result, we can introduce the measure $\cL$ on the skeleton of $\cT$ in the same way as in \eqref{def: L}, with $\beta$ as the quadratic coefficient in $\Psi$. Pruning $\cT$ using a Poisson point measure of intensity $\cL$ is considered in \cite{AbDe08,AbDe10}. Moreover, if we denote by $\cT_{t}$ the connected component containing the root at time $t$, then it is shown in \cite{AbDe12} that $\cT_{t}$ under $\mathbf N^{\Psi}$ is distributed as a $\Psi_{t}$-L\'evy tree, where $\Psi_{t}(\lambda)=\Psi(\lambda+t)-\Psi(t)$.   It follows that under $\mathbf N^{\Psi}$, the law of $\cT_{t}$ is absolutely continuous to that of $\cT$. 
The cut tree associated to this pruning of $\cT$ can be defined in the same way as for the ICRT. Write $\cC$ for the cut tree. It turns out that we once again has a distributional identity: $\cC\eqd \cT$ under $\mathbf N^{\Psi}$; this result is essentially due to Abraham and Delmas \cite{AbDe13}. 

To apply the previous strategy in reconstructing $\cT$ from $\cC$, one will need three ingredients. First, an analogue to the distance approximation in \eqref{eq: icrt-app}; this, in fact, is known as a direct consequence of Bismut's decomposition for L\'evy trees \cite{DuLeG05}. Secondly, a result on the local time preservation in the spirit of Proposition~\ref{prop: br_pt}. Thirdly, the definition and distribution of the collection of traces. In view of the results in \cite{AbDeHo14}, it is reasonable to expect that these two missing ingredients are not far from reach.

\section{Further properties of inhomogeneous continuum random trees}

In this section, we prove the approximations for the local time of a branch point in 
Proposition~\ref{prop:degrees}, and for 
the distance between two points in 
Proposition~\ref{prop: d-app}. Both rely on the Line-breaking construction recalled in Section~\ref{sec:line_breaking}. 

\subsection{Local time of a branch point: Proof of Proposition~\ref{prop:degrees}}
\label{sec: br-pt}

Recall the real tree $\cR_{k}$ with the first $k$ branches in the Line-breaking construction in Section~\ref{sec:line_breaking}. 
We denote by $\cR_{\infty}=\cup_{k\ge 1}\cR_{k}\subseteq \cT$. It is not difficult to check that $\Br(\cR_{\infty})=\Br(\cT)$ and therefore $\Br_{\infty}(\cR_{\infty})=\Br_{\infty}(\cT)$ a.s.
On the other hand, it is clear from the Line-breaking Algorithm that $\Br_{\infty}(\cR_{\infty})=\{B_{i}: \theta_{i}>0\}$, where we recall that $B_{i}$ is the image of $\xi_{i, 1}$ in $\cT$. 
This proves the first statement in Proposition~\ref{prop:degrees}. 
For the second, we observe that if $\sigma\in \cT$ has bounded degree, then $\deg(\sigma, \cR_{k})/\Psi^{-1}(k)\to 0$ a.s. 
It remains to show that for each $\theta_{i}>0$,
\begin{equation}
\label{def: br-pt'}
\frac{\deg(B_{i}, \cR_{k})}{\Psi^{-1}(k)} \to \theta_{i}, \quad \text{in probability}.
\end{equation}
To that end, recall that $(U_{j})_{j\ge 1}$ is a Poisson process of rate $\beta x\,dx$. 
Recall also that $(\xi_{i, j})_{j\ge 1}$ is a Poisson process of rate $\theta_{i}$. 
Let us also define for $t\ge 0$, 
\[
N_{i}(t)=\sum_{j\ge 1}\mathbf 1_{\{\xi_{i, j}\le t\}}=\max\{j\ge 1: \xi_{i, j}\le t\}. 
\]
Recall that $\cR_{k}$ is obtained from $\cR_{k-1}$ by identifying $\eta_{k-1}$ with $\eta_{k-1}^{\ast}$. 
It follows that  
\[
\deg(B_{i}, \cR_{k}) = N_{i}(\eta_{k-1})+1, \quad \text{if } N_{i}(\eta_{k-1})\ge 2. 
\] 
On the other hand, Law of Large Numbers entails that $N_{i}(t)/t\to \theta_{i}$ almost surely. 
Together with the fact  $\eta_{k}\to\infty$ a.s., we deduce that
\[ 
\frac{\deg(B_{i}, \cR_{k})}{\eta_{k-1}}\xrightarrow{k\to\infty} \theta_{i}, \quad \text{almost surely}.
\]
Therefore, to obtain \eqref{def: br-pt'}, it suffices to show that 
\begin{equation}
\label{eta-asy}
\frac{\eta_{k-1}}{\Psi^{-1}(k)}\xrightarrow{k\to\infty} 1 \quad \text{in probability.}
\end{equation}
Its proof will occupy the rest of this part. Let us first introduce:  
\[
\sL(t)=\sum_{k\ge 1}\mathbf 1_{\{\eta_{k}\le t\}}=\sum_{i\ge 1}(N_{i}(t)-1)_{+}+\sum_{j\ge 1}\mathbf 1_{\{U_{j}\le t\}} \quad \text{and}\quad  a(t)=\beta t+\sum_{i\ge 1}\theta_{i}\mathbf 1_{\{\xi_{i, 1}\le t\}}. 
\]
Properties of Poisson processes imply that 
\[
M^{0}_{t}:=\sum_{j\ge 1}\mathbf 1_{\{U_{j}\le t\}}-\int_{0}^{t}\beta s ds, \quad t\ge 0,
\] 
is a martingale with $\mathbb E[(M^{0}_{t})^{2}]=\int_{0}^{t}\beta sds$. On the other hand, for each $i\ge 1$, 
\[
M^{i}_{t}:=(N_{i}(t)-1)_{+}-\int_{0}^{t}\theta_{i}\mathbf 1_{\{\xi_{i, 1}\le s\}}ds =(N_{i}(t)-1)_{+}-\theta_{i}(t-\xi_{i, 1})_{+} , \quad t\ge 0,
\]
is also a martingale, whose second moment can be bounded as follows:
\[
\mathbb E[(M^{i}_{t})^{2}] \le  \mathbb E[(N_{i}(t)-1)_{+}^{2}] = (\theta_{i}t)^{2}-(e^{-\theta_{i}t}+\theta_{i}t-1)\le (\theta_{i}t)^{2}.
\]
Since $M^{i}_{t}, i\ge 0$, are independent, it follows that
\[
\frac{1}{\Psi(t)^{2}}\mathbb E\Big[\Big(\sL(t)-\int_{0}^{t}a(s)ds\Big)^{2}\Big]=\frac{1}{\Psi(t)^{2}}\sum_{i\ge 0}\mathbb E\big[(M^{i}_{t})^{2}\big]\le \frac{1}{\Psi(t)^{2}}(\|\bth\|^{2}+\beta)t^{2}. 
\]
We observe that the above tends to 0 as $t\to\infty$, 
since $\Psi(t)/t$ diverges thanks to \eqref{Theta-ass}. 
As a consequence, 
\[
\frac{1}{\Psi(t)}\Big(\sL(t)-\int_{0}^{t}a(s)ds\Big) \xrightarrow{t\to\infty} 0, \quad \text{in probability}.
\]
On the other hand, noting  that $\Psi'(t)=\mathbb E[a(t)]$, we have
\[
\mathbb E\bigg[\bigg(\frac{a(t)}{\Psi'(t)}-1\bigg)^{2}\bigg]=\frac{\mathrm{Var}(a(t))}{\Psi'(t)^{2}}=\frac{1}{\Psi'(t)^{2}} \sum_{i\ge 1}\theta_{i}^{2}(1-e^{-\theta_{i}t})e^{-\theta_{i}t} \le \left(\frac{\|\bth\|}{\Psi'(t)}\right)^{2}\xrightarrow{t\to\infty} 0. 
\]
Hence, $a(t)/\Psi'(t)\to 1$ in probability. By L'H\^opital's rule, this implies that
\[
\frac{1}{\Psi(t)}\int_{0}^{t}a(s)ds \xrightarrow{t\to\infty} 1, \quad \text{in probability}. 
\]
As a consequence, we obtain that 
\begin{equation}
\label{eq: Lt}
\frac{\sL(t)}{\Psi(t)}=\frac{1}{\Psi(t)} \Big(\sL(t)-\int_{0}^{t}a(s)ds\Big)+\frac{1}{\Psi(t)}\int_{0}^{t}a(s)ds \xrightarrow{t\to\infty} 1, \quad \text{in probability.}
\end{equation}
Now let us show how this implies \eqref{eta-asy}. For $\epsilon>0$, let us denote $t_{k}=\Psi^{-1}((1+\epsilon)k)$. Let $\delta\in (0, \epsilon/(1+\epsilon))$, so that $k<(1-\delta)\Psi(t_{k})$. It follows that 
\[
\mathbb P(\eta_{k}>t_{k})=\mathbb P(\sL(t_{k})<k)\le \mathbb P\big(\sL(t_{k})\le (1-\delta)\Psi(t_{k})\big)\to 0, \quad k\to\infty,
\]
as a consequence of \eqref{eq: Lt}. Denoting $t_{k}'=\Psi^{-1}((1-\epsilon)k)$, we can show similarly that $\mathbb P(\eta_{k}\le t_{k}')\to 0$. As $\Psi$ is convex and increasing, $\Psi^{-1}$ is concave and increasing, with $\Psi^{-1}(0)=0$. This leads to the following bound: for all $\epsilon>0$ and $t>0$, 
\begin{equation}
\label{bd: psiin}
\Psi^{-1}\big((1+\epsilon)t\big)\ge \Psi^{-1}(t)=\Psi^{-1}\Big(\frac{1}{1+\epsilon}(1+\epsilon)t + \big(1-\frac{1}{1+\epsilon}\big)\cdot 0\Big)\ge \frac{1}{1+\epsilon}\Psi^{-1}\big((1+\epsilon)t\big).
\end{equation}
We apply this to $t_{k}$ and $t_{k}'$ and find that for all $k\ge 1$, 
\[
t_{k}=\Psi^{-1}\big((1+\epsilon)k\big)\le (1+\epsilon) \Psi^{-1}(k), \quad t_{k}'=\Psi^{-1}\big((1-\epsilon)k\big) \ge (1-\epsilon)\Psi^{-1}(k).
\]
Together with previous arguments, this shows that $\mathbb P((1-\epsilon)\Psi^{-1}(k)\le \eta_{k}\le (1+\epsilon)\Psi^{-1}(k))\to 1$, as $k\to\infty$. Using \eqref{bd: psiin} again, we see that $\Psi^{-1}(k)\le \Psi^{-1}(k+1)\le \frac{k+1}{k}\Psi^{-1}(k)$, so that $\Psi^{-1}(k+1)/\Psi^{-1}(k)\to 1$. 
This completes the proof of \eqref{def: br-pt'}. Finally, we deduce from \eqref{id: re-root'} that
\[
\big(\cR'_{k}, \deg(B_{i}, \cR'_{k})\big)  \eqd \big(\cR_{k}, \deg(B_{i}, \cR_{k})\big). 
\]
As $\cR'_{k-1}$ can be recovered from $\cR'_{k}$ by removing a uniform leaf, the above identity leads to the following: for all $m\in \N$ and integers $k_{1}\ge k_{2}\ge \dots k_{m}$, we have 
\[
\big(\deg(B_{i}, \cR'_{k_{j}})\big)_{1\le j\le m} \eqd  \big(\deg(B_{i}, \cR_{k_{j}})\big)_{1\le j\le m}.
\]
This implies  $(\deg(B_{i}, \cR'_{k}))_{k\ge 1}$ has the same law as $(\deg(B_{i}, \cR_{k}))_{k\ge 1}$. The proof of Proposition~\ref{prop:degrees} is now complete. 

\subsection{Length of a uniform branch: Proof of Proposition~\ref{prop: d-app}}
\label{sec: icrt-app}


Let us recall some notation introduced in Section~\ref{sec:line_breaking} about the Line-breaking construction: $(\xi_{i, j})_{j\ge 1}$ is a Poisson process of rate $\theta_{i}$, $i\ge 1$. In particular, $(\xi_{i, 1})_{i\ge 1}$ is a collection of independent exponential variables with $\mathbb E[\xi_{i, 1}]=\theta_{i}^{-1}$. 
Thanks to \eqref{id: re-root'}, we observe that
the convergence in \eqref{eq: icrt-app}  reduces to the following
\begin{equation}
\label{eq: icrt-app'}
\frac{1}{\gamma_{\bth}(\epsilon)}\sum_{i: \theta_{i}>\epsilon}\mathbf 1_{\{\xi_{i, 1}\le \eta_{1}\}} \xrightarrow{\epsilon\to 0} \eta_{1}, \quad \text{in probability}.
\end{equation}
Let us show this. To ease the notation, 
from now on we will write $\gamma(\epsilon)$ instead of $\gamma_{\bth}(\epsilon)$. 
Let us introduce the following quantities: for $t\ge 0$ and $\epsilon>0$, 
\[
L_{\epsilon}(t)= \sum_{i: \theta_{i}>\epsilon} \mathbf 1_{\{\xi_{i, 1}\le t\}} \quad \text{and, for each }i\ge 1, \quad  M_{i}(t)=\mathbf 1_{\{\xi_{i, 1}\le t\}}-\theta_{i}(t\wedge \xi_{i, 1})\,. 
\]
It is straightforward to check that, for each $i\ge 1$, $(M_{i}(t))_{t\ge 0}$ is a martingale with respect to the natural filtration of $(\mathbf 1_{\{\xi_{i, 1}\le t\}})_{t\ge 0}$ and that $\mathbb E[M_{i}(t)^{2}]\le \mathbb P(\xi_{i, 1}\le t)=1-\exp(-\theta_{i}t)$. 
Now let
\[
A_{\epsilon}(t)=\sum_{i: \theta_{i}>\epsilon}\theta_{i}(\xi_{i, 1}\wedge t) \quad \text{and}\quad \mathcal M_{\epsilon}(t)=\sum_{i: \theta_{i}>\epsilon}M_{i}(t), \quad t\ge 0. 
\]
Then $(\mathcal M_{\epsilon}(t))_{t\ge 0}$ is a martingale with respect to the natural filtration of $\{\mathbf 1_{\{\xi_{i, 1}\le t\}}, t\ge 0, i\ge 1\}$. 
Doob's maximal inequality applies, so that
\begin{equation}
\label{bd: mart}
\frac{1}{\gamma(\epsilon)^{2}}\mathbb E\Big[\sup_{s\le t}\mathcal M_{\epsilon}(s)^{2}\Big]\le \frac{4}{\gamma(\epsilon)^{2}}\mathbb E\Big[\mathcal M_{\epsilon}(t)^{2}\Big]\le \frac{4}{\gamma(\epsilon)^{2}}\sum_{i: \theta_{i}>\epsilon}(1-e^{-\theta_{i}t}) \le \frac{4t}{\gamma(\epsilon)},
\end{equation}
where in the equality above we have used the independence of $M_{i}(t)$, $i\ge 1$. On the other hand, we observe that 
\[
\bigg|\frac{1}{\gamma(\epsilon)}A_{\epsilon}(t)-t\bigg|
=\bigg|\frac{1}{\gamma(\epsilon)}\sum_{i: \theta_{i}>\epsilon}\theta_{i}\int_{0}^{t}\big(\mathbf 1_{\{\xi_{i, 1}> s\}}-1\big)ds\bigg|
=\frac{1}{\gamma(\epsilon)}\sum_{i: \theta_{i}>\epsilon}\int_{0}^{t}\theta_{i}\mathbf 1_{\{\xi_{i, 1}\le s\}}ds,
\]
which is clearly increasing in $t$. It follows that
\begin{align}
\notag
\mathbb E\bigg[\sup_{s\le t}\bigg|\frac{1}{\gamma(\epsilon)}A_{\epsilon}(s)-s\bigg|\bigg] &= \mathbb E\bigg[\frac{1}{\gamma(\epsilon)}\sum_{i: \theta_{i}>\epsilon}\int_{0}^{t}\theta_{i}\mathbf 1_{\{\xi_{i, 1}\le s\}}ds\bigg] \\ \label{bd: comp}
&\le\frac{1}{\gamma(\epsilon)}\sum_{i: \theta_{i}>\epsilon}\theta_{i}\int_{0}^{t}(1-e^{-\theta_{i}s})ds\le \frac{\|\bth\|^{2}t^{2}}{\gamma(\epsilon)}. 
\end{align}
Note that $L_{\epsilon}(t)=A_{\epsilon}(t)+\mathcal M_{\epsilon}(t)$. Therefore, \eqref{bd: mart} and \eqref{bd: comp} yield the following:  
\begin{equation}\label{cv: L} 
\mathbb E\bigg[\sup_{s\le t}\bigg|\frac{L_{\epsilon}(s)}{\gamma(\epsilon)}-s\bigg|\bigg] \le \frac{\|\bth\|^{2}t^{2}}{\gamma(\epsilon)}+\frac{2\sqrt{t}}{\sqrt{\gamma(\epsilon)}}\to 0,  
\end{equation}
as $\epsilon\to 0$, since $\gamma(\epsilon)\to \infty$ by assumption. 
On the other hand, the law of $\eta_{1}$ is tight: $\mathbb P(\eta_{1}>t)\to 0$ as $t\to\infty$. Combined with \eqref{cv: L},
this allows us to deduce \eqref{eq: icrt-app'}.

Now suppose that $\theta_{i}, i\ge 1$ are all distinct. Let us show that the convergence in \eqref{eq: icrt-app'} also takes place almost surely. Note that in this case there is at most one $\theta_{i}=\epsilon$; therefore $ \gamma(\epsilon)-\gamma(\epsilon+)\le \epsilon$ for all $\epsilon>0$. Since  $\gamma(\epsilon)\to \infty$ as $\epsilon\to 0$, it is possible to choose a sequence $\epsilon_{m}\to 0$ so that for $m$ sufficiently large, 
\begin{equation}
\label{bd: subseq}
m^{\delta}\le \gamma(\epsilon_{m})<(m+1)^{\delta}, 
\end{equation}
where $\delta\ge 1$ is some real number to be specified later. 
Indeed, let $\epsilon_{m}=\inf\{\epsilon: \gamma(\epsilon)<m^{\delta}\}$; then by definition, $\gamma(\epsilon_{m}+)\le m^{\delta}$, so that $\gamma(\epsilon_{m})\le m^{\delta}+\epsilon_{m}\le (m+1)^{\delta}$ for $m$ sufficiently large. 
Meanwhile, we note that the function 
\[
g(t):=\tfrac12\beta t^{2}+\sum_{i\ge 1}\big(\theta_{i}t-\log(1+\theta_{i}t)\big)
\]
is convex on $(0, \infty)$ since $g''(t)=\beta+\sum_{i\ge 1}\frac{\theta_{i}^{2}}{(1+\theta_{i}t)^{2}}\in (0, \infty)$. As $g(0)=0$, Jensen's inequality yields that $g(1)=g(\tfrac1t\cdot t+(1-\tfrac1t)\cdot0)\le \tfrac1t g(t)$, for all $t\ge 1$. 
Recall from the Line-breaking construction that $\eta_{1}=U_{1}\wedge \min_{i\ge 1}\xi_{i, 2}$. We deduce that
\begin{equation}
\label{bd: eta}
\mathbb P(\eta_{1}>t)=\exp(-g(t))\le \exp(-Ct), \quad \text{ for all } t\ge 1,
\end{equation}
where we have denoted $C=C(\beta, \bth)=g(1)\in (0, \infty)$. 
Since 
\[
\mathbb P\bigg(\,\bigg|\frac{L_{\epsilon_{m}}(\eta_{1})}{\gamma(\epsilon_{m})}-\eta_{1}\bigg| > \frac{1}{m}\bigg)\le \mathbb P(\eta_{1}>m) + \mathbb P\bigg(\sup_{s\le m}\bigg|\frac{L_{\epsilon_{m}}(s)}{\gamma(\epsilon_{m})}-s\bigg| > \frac{1}{m}\bigg), 
\]
bounding the first term using \eqref{bd: eta}, while applying  \eqref{cv: L}, \eqref{bd: subseq} and the Markov inequality for the second, we find that 
\[
\mathbb P\bigg(\,\bigg|\frac{L_{\epsilon_{m}}(\eta_{1})}{\gamma(\epsilon_{m})}-\eta_{1}\bigg| > \frac{1}{m}\bigg)
\le e^{-Cm}+\|\bth\|^{2}m^{3-\delta}+2m^{\frac{3-\delta}{2}}.
\]
Now take $\delta\ge 7$.  The Borel--Cantelli Lemma is in force. Hence, 
\[
\frac{L_{\epsilon_{m}}(\eta_{1})}{\gamma(\epsilon_{m})} \xrightarrow{m\to\infty} \eta_{1}, \quad \mathbb P\text{-a.s.}
\]
For $\epsilon\in [\epsilon_{m+1}, \epsilon_{m})$, the monotonicity of both $L_{\epsilon}(t)$ and $\gamma(\epsilon)$ yields
\[
\frac{L_{\epsilon_{m}}(\eta_{1})}{\gamma(\epsilon_{m})}\frac{\gamma(\epsilon_{m})}{\gamma(\epsilon_{m+1})}=\frac{L_{\epsilon_{m}}(\eta_{1})}{\gamma(\epsilon_{m+1})} \le \frac{L_{\epsilon}(\eta_{1})}{\gamma(\epsilon)}\le \frac{L_{\epsilon_{m+1}}(\eta_{1})}{\gamma(\epsilon_{m})}=\frac{L_{\epsilon_{m+1}}(\eta_{1})}{\gamma(\epsilon_{m+1})}\frac{\gamma(\epsilon_{m+1})}{\gamma(\epsilon_{m})}. 
\]
By \eqref{bd: subseq}, $1\le \gamma(\epsilon_{m+1})/\gamma(\epsilon_{m})\le (m+2)^{\delta}/m^{\delta}\to 1$, as $m\to \infty$. Therefore, both sides above converge to $\eta_{1}$ almost surely, in consequence of the previous convergence. 
This completes the proof of Proposition~\ref{prop: d-app}. 


\section{Further properties of cut trees: Proof of Proposition~\ref{prop: br_pt}}

\label{sec:proof_correspondence}

In this section, we prove Proposition~\ref{prop: br_pt}, which shows that the map $\phi: \Br(\cC)\to \cP$ preserves the local times. The arguments below rely on the approximation of local times as found in Proposition~\ref{prop:degrees}. 

Recall that $(V_{i})_{i\ge 0}$ is a sequence of i.i.d.~points in $\cT$ with common law $\mu$. Let $\cR'_{k}=\cup_{1\le i\le k}\llb V_{0}, V_{i}\rrb$, which has the same distribution as the tree $\cR_{k}$ that appear in the Line-breaking construction. Recall also the notation $\cS_{k}=\cup_{1\le i\le k}\llb 0, i\rrb$ for the spanning tree in $\cC$. 
Fix $\sigma\in \Br(\cC)$. Let us write $D_{\sigma}(k)=\deg(\sigma, \cS_{k})$ and $\tilde D_{\sigma}(k)=\deg(\phi(\sigma), \cR_{k})$, with the convention that if $\sigma\notin \cS_{k}$, then $D_{\sigma}(k)=-\infty$ (and similarly for $\tilde D_{\sigma}(k)$). 
In view of  Proposition~\ref{prop:degrees}, we aim to compare $D_{\sigma}(k)$ with $\tilde D_{\sigma}(k)$. 

Let us recall the equivalence relation $\sim_{t}$ on $[k]=\{1, 2, \dots, k\}$ induced by removing the points of $\cP_{t}$ from $\cR'_{k}$. Recall also that $\Pi^{[k]}_{t}$ is the partition of $[k]$ associated to $\sim_{t}$. Let us write $\#\Pi^{[k]}_{t}$ for the number of equivalence classes in this partition. We note that there is a unique time $\tau(\sigma)$ so that $(\tau(\sigma), \phi(\sigma))\in \cP$. The effect of removing $\phi(\sigma)$ on the partition-valued process $(\Pi^{[k]}_{t})_{t\ge 0}$ can be easily understood: there is an equivalence class, say $\pi\in \Pi^{[k]}_{\tau(\sigma)-}$, which disappears at time $\tau(\sigma)$ and is replaced by a number of new  equivalence classes. Denote by $Q_{k}(\sigma)$ the number of new equivalence classes. Note that $Q_{k}(\sigma)$ corresponds to the number of children of $\sigma$ in $\cS_{k}$. This yields the following:
\begin{equation}
\label{bd: deg}
Q_{k}(\sigma)+1 = D_{\sigma}(k) \quad \text{and}\quad Q_{k}(\sigma) - 1 = \#\Pi^{[k]}_{\tau(\sigma)}-\#\Pi^{[k]}_{\tau(\sigma)-}.
\end{equation}
On the other hand, let $\cR'(\sigma, k)$ be the connected component of $\cR'_{k}\setminus\cP_{\tau(\sigma)-}$ that contains $\phi(\sigma)$; then it is not difficult to see from the pruning that $Q_{k}(\sigma)$ corresponds to the number of connected components of $\cR'(\sigma, k)\setminus\{\phi(\sigma)\}$ that contains some $V_{i}$, $1\le i\le k$. In particular, $Q_{k}(\sigma)$ is bounded by the degree of $\phi(\sigma)$ in $\cR'_{k}$. Combined with \eqref{bd: deg}, this leads to 
\begin{equation}
\label{bd: deg'}
D_{\sigma}(k)\le \tilde D_{\sigma}(k) + 1. 
\end{equation}
Observe that this is the counterpart of \eqref{ineq-degree}. 

\paragraph{Case $\sigma\in \Br_{\infty}(\cC)$.}  In this case, $D_{\sigma}(k)\nearrow \infty$. Dividing both sides of the inequality in \eqref{bd: deg'} by $\Psi^{-1}(k)$ and letting $k\to\infty$, we find that
\begin{equation}
\label{bd: deg2}
\Delta_{\cC}(\sigma) \le \Delta_{\cT}\big(\phi(\sigma)\big).
\end{equation}
Recall that $\cC$ is distributed as an ICRT with parameter $(\beta, \bth)$. Rank the infinite-degree branch points in $\cC$ by their local times and denote by $\tilde B_{i}$ the one with $\Delta_{\cC}(\tilde B_{i})=\theta_{i}$. 
Plugging $\sigma=\tilde B_{1}$ into \eqref{bd: deg2} yields $\Delta_{\cT}(\phi(\tilde B_{1}))=\theta_{1}$. Repeating this procedure with $\tilde B_{2}, \tilde B_{3} \dots$ and bearing in mind that $\phi$ is injective from its definition, we see that 
\[
\Delta_{\cC}(\sigma)=\Delta_{\cT}\big(\phi(\sigma)\big) \quad \text{for all } \sigma\in \Br_{\infty}(\cC).
\]

\paragraph{Case $\sigma\notin \Br_{\infty}(\cC)$. } Let us show in this case that $\deg(\phi(\sigma), \cT)<\infty$. Since $\cP_{\infty}$ is disjoint from the leaf set or the set of binary branch points of $\cT$, this will allow us to conclude that $\deg(\phi(\sigma), \cT)=2$ a.s. 
Write $\mathcal A$ for the event $\{\phi(\sigma)\in \Br_{\infty}(\cT)\}$. 
On this event, let $\cT^{1}, \cT^{2}, \dots $ be the connected components of $\cT\setminus\{\phi(\sigma)\}$ ranked in decreasing order of their $\mu$-masses. 
Note that $(V_{i})_{i\ge 1}$ is dense almost everywhere in $\cT$, since the support of $\mu$ is $\cT$. 
This means that we can choose a subsequence $(V_{i_{k}})_{k\ge 1}$ of $(V_{i})$ satisfying $V_{i_{k}}\in \cT^{k}$ and $V_{i_{k}}\to \phi(\sigma)$ a.s. 
If $\deg(\sigma, \cC)$ is bounded, then $(Q_{k}(\sigma))_{k\ge 1}$ is eventually constant, and therefore, apart from a finite number the remaining $V_{i_{k}}$ are all disconnected from $\phi(\sigma)$ before $\tau(\sigma)$. 
It follows that for all $t\ge 0$, 
\begin{align}
\label{pf: comp}
\mathbb P\big(\deg(\sigma, \cC)<\infty\,|\,\mathcal A\big) \le \mathbb P\big(\tau(\sigma)>t\,|\, \mathcal A\big)
+ \lim_{m\to\infty}\mathbb P\big(\llb V_{i_{k}}, \phi(\sigma)\llb \, \cap \,\cP_{t}\ne \varnothing, \forall \, k\ge m\,|\,\mathcal A \big)
\end{align}
By first conditioning on $\cT$, we can re-write the second term as follows:
\[
\mathbb P\big(\llb V_{i_{k}}, \phi(\sigma)\llb \, \cap \cP_{t}\ne \varnothing, \forall \, k\ge m\,|\,\mathcal A \big) =\mathbb E\Big[\prod_{k\ge m}\Big(1-e^{-t \cL(\llb V_{i_{k}}, \phi(\sigma)\llb)}\Big)\,\Big|\,\mathcal A \Big].
\]
However, since $d(\phi(\sigma), V_{i_{k}})\to 0$ a.s., we also have $\cL(\llb V_{i_{k}}, \phi(\sigma)\llb) \to 0$ a.s. 
In particular, we can find some $M\in (0, \infty)$ so that $\mathbb P(\sup_{k}\cL(\llb V_{i_{k}}, \phi(\sigma)\llb)>M)\le \epsilon$. 
Therefore,
\[
\mathbb P\big(\llb V_{i_{k}}, \phi(\sigma)\llb \, \cap\, \cP_{t}\ne \varnothing, \forall \, k\ge m\,|\,\mathcal A \big)  \le \epsilon + (1-e^{-tM})^{\infty} = \epsilon. 
\]
For those $\phi(\sigma)\in \Br_{\infty}(\cT)$, the law of $\tau(\sigma)$ is tight. Therefore, by first letting $\epsilon\to0$ and then letting $t\to\infty$ in \eqref{pf: comp}, we conclude that $\mathbb P\big(\deg(\sigma, \cC)<\infty\,|\,\phi(\sigma)\in \Br_{\infty}(\cT)\big)=0$. Since we have assumed $\sigma\notin \Br_{\infty}(\cC)$, it can only happen that $\mathbb P(\phi(\sigma)\in\Br_{\infty}(\cT))=0$.

\section{Proof of Proposition~\ref{prop: route}}
\label{sec: trace-pf}

\subsection{Strategy of the proof}

In this part, we give a proof of Proposition~\ref{prop: route} by constructing a collection of traces that meets the requirement. 
To do so, we will rely upon the following observation. 
Let $(E, d)$ be a {\it complete} metric space and suppose that $(u_{i})_{i\ge 1}$ is dense and pairwise distinct in $E$. 
Then any point $x\in E$ is characterised by the sequence $\delta(x)=(\delta_{i})_{i\ge 1}$, where we have put $\delta_{i}=d(x, u_{i})$. Observe that the sequence $\delta(x)=(\delta_{i})_{i\ge 1}$ has the following properties:
\begin{compactenum}
\item[(A1)]
for all $i\ge 1$, $\delta_{i}\in [0, \infty)$;
\item[(A2)]
$\inf_{i\ge 1}\delta_{i}=0$;
\item[(A3)]
for $i, j\ge 1$, $|\delta_{i}-\delta_{j}|\le d(u_{i}, u_{j})\le \delta_{i}+\delta_{j}$. 
\end{compactenum}

\smallskip
As it turns out, the converse is also true. 
\begin{lem}
\label{lem: char}
If $\delta=(\delta_{i})_{i\ge 1}$ is a sequence that satisfies (A1-A3), then there is a unique point $y\in E$ so that $\delta_{i}=d(y, u_{i})$, $i\ge 1$. 
\end{lem}

\begin{proof}
According to (A2), we can find a sequence $(i_{n})_{n\ge 1}$ so that $\delta_{i_{n}}<2^{-n}$, $n\ge 1$. Applying (A3) to the sequence yields that $d(u_{i_{n}}, u_{i_{m}})\le 2^{-n+1}$, for all $m\ge n$. Therefore, $(u_{i_{n}})_{n\ge 1}$ is a Cauchy sequence, which converges to some point $y\in E$ as $E$ is complete. 
Appealing to (A3) once again, we find $|\delta_{i}-\delta_{i_{n}}|\le d(u_{i}, u_{i_{n}})\le \delta_{i}+\delta_{i_{n}}$. Letting $n\to\infty$ leads to $\delta_{i}= d(u_{i}, y)$, $i\ge 1$. For the uniqueness, say $y'\in E$ also satisfies $\delta_{i}=d(y', u_{i})$, $i\ge 1$. We can then find a sequence $i'_{n}\to \infty$ so that $u_{i'_{n}}\to y'$ and $\delta_{i'_{n}}\to 0$. By (A3), $d(u_{i'_{n}}, u_{i_{n}})\le \delta_{i'_{n}}+\delta_{i_{n}}$. Taking $n\to\infty$ concludes the proof. 
\end{proof}

\noindent{\bf Remark. }
Equipping $\R^{\infty}$ with the product topology, we see from the previous proof that the map from $\delta=(\delta_{i})_{i\ge 1}\in \R^{\infty}$ to the unique point $y\in E$ that satisfies $\delta_{i}=d(y, u_{i})$ is measurable. 

\medskip

Let $\mathscr P(E)$ stand for the set of (Borel) probability measures on $(E, d)$, equipped with the weak topology. 
Similarly, let $\mathscr P(\R^{\infty})$ be the set of probability measures of $\R^{\infty}$ with respect to the product topology on $\R^{\infty}$. 
The following lemma will be used implicitly later on.

\begin{lem}
Define a map $f: \mathscr P(E)\to \mathscr P(\R^{\infty})$ as follows: for $\mu\in \mathscr P(E)$, $f(\mu)$ is the law of $ (d(u_{i}, U))_{i\ge 1}$, where $U$ is a random point with law $\mu$. Then $f$ is injective and continuous. 
\end{lem}

\begin{proof}
Injectivity is clear from Lemma~\ref{lem: char}. For the continuity, consider a sequence of random variables $(U_{n})_{n\in\N}$ converging to some limit $U$ in distribution. Then we have $(d(u_{i}, U_{n}))_{1\le i\le k}\to (d(u_{i}, U))_{1\le i\le k}$ in distribution for any $k$. This proves the statement. 
\end{proof}

\smallskip
\noindent
{\bf An overview of the proof for Propositions~\ref{prop: route} and~\ref{prop: tr-dist}.} 
Rather than constructing the desired collection $\mathcal W$ directly, we will first introduce a collection $\mathcal R(i, j)$ which is the candidate for  $\mathcal R_{\mathcal W}(i, j)$,  $(i, j)\in \N^{2}$. We recall that the collection $\mathcal R_{\mathcal W}(i, j)$ is indexed by the internal points of a complete binary tree $\mathbb U^{\ast}=\cup_{n\ge 1}\{0, 1\}^{n}$. This reflects the fact that every cut on the path connecting $V_{i}$ to $V_{j}$ in $\cT$ splits it into two further segments (see Fig~\ref{fig4} below for an illustration). To track the fragmentation of $\llb V_{i}, V_{j}\rrb$ into disjoint segments, we introduce a separate collection $\mathcal Y(i, j) = (y^{i, j}_{u})_{u\in \mathbb U^{\ast}}$ formed by the endpoints of these segments. The idea is that each element $r^{i, j}_{u}\in \mathcal R(i, j)$ will act as the ``trace'' of $y^{i, j}_{u}\in \cT$ in a similar way that $i\in \cC$ is the ``trace'' of $V_{i}\in \cT$. Drawing a parallel with the definition \eqref{def: delta'} of $\delta(i, j)$, it is then natural to expect that the distance between $r^{i, j}_{u}$ and $k\in \N$, denoted as $\delta^{i, j}_{u}(k)$, is given as follows:
\begin{equation}
\label{dist: delta''}
\delta^{i, j}_{u}(k) = \int_{\tau^{i, j}_{u, k}}^{\infty} \Big\{\mu^{i, j}_{t}(u)+\mu\big(\cT_{t}(V_{k})\big)\Big\}dt,
\end{equation}
where $\tau^{i, j}_{u, k}$ is the moment when $y^{i, j}_{u}$ is disconnected from $V_{k}$, and $\mu^{i, j}_{t}(u)$ represents the $\mu$-mass of the part of $\cT$ containing $y^{i, j}_{u}$. However, defining the latter requires some more work, as $y^{i, j}_{u}$ could be the location of a cut; in that case, $y^{i, j}_{u}$ is removed from the tree at the moment of the cut. The removal of $y^{i, j}_{u}$ splits $\cT$ into some further connected components, and it turns out that for each of these components, there is a natural ``copy'' of $y^{i, j}_{u}$, and tracking the mass of the component containing this copy will allow us to define properly the quantity $\mu^{i, j}_{t}(u)$. 

Thanks to Lemma~\ref{lem: char}, once we have checked that $(\delta^{i, j}_{u}(k))_{k\ge 1}$ verifies (A1-A3), we will be able to identify the point $r^{i, j}_{u}$ in $\cC$. This will be done using weak convergence arguments, by looking at the counterparts of $y^{i, j}_{u}$ and $r^{i, j}_{u}$ in the $\bp$-trees.  
Following the discussion in Section~\ref{sec: intro-p-tree}, it should become clear that in a $\bp$-tree, the role of $y^{i, j}_{u}$ is played by a neighbour to the vertex where a cut falls. Since the cut tree shares the same vertex set as the $\bp$-tree, the role of $r^{i, j}_{u}$ is then played by the same vertex in the cut tree. 
It is also possible to determine the distribution of the discrete counterpart of $y^{i, j}_{u}$, thanks to the product form of $\pi_{\bp}$ in \eqref{def: pi}. From there we deduce the distribution of the discrete counterpart of $r^{i, j}_{u}$. 
This analysis is actually done in \cite{BW}; we will recall the main conclusions in Lemma~\ref{lem: p-tree}. 
We then write down an analogue of \eqref{dist: delta''}, which corresponds to the distance between the discrete counterpart of $r^{i, j}_{u}$ and a sequence of i.i.d.~vertices in the cut tree. 
Finally, we show that jointly with the convergence of the pruning processes on $\bp$-trees to the pruning of ICRT, the discrete counterpart of \eqref{dist: delta''} converges to its continuum analogue, which must have the distribution of the pairwise distances between certain random point and the elements of $\N$ in $\cC$. 
Note that this approach will also allow us to determine the distribution of $r^{i, j}_{u}$. 

Having constructed the routings $(\mathcal R(i, j))_{i, j\in \N}$, we will proceed to show that they satisfy certain consistency conditions, so that we will be able to extract the collection $\mathcal W$ from the routings. 

\subsection{Constructing the routing $\mathcal R(i, j)$: first generation}

\begin{figure}
\centering
\includegraphics[height = 5cm]{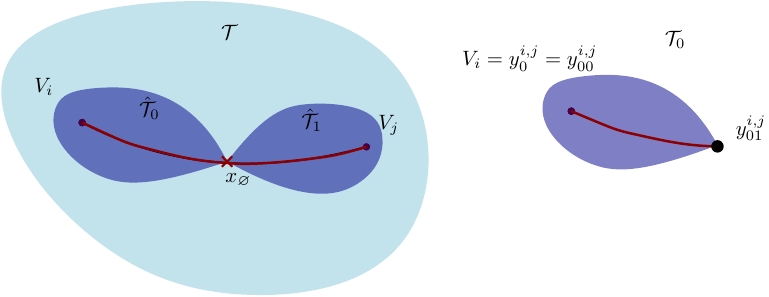}
\caption{\label{fig4} On the left, the red curve depicts the path between $V_{i}$ and $V_{j}$ in $\cT$. The cross represents the location of the first cut $x_{\varnothing}$ on the path. The darkly shaded areas correspond respectively to $\hat{\mathcal T}_{0}$ and $\hat{\mathcal T}_{1}$, the connected components containing respectively $V_{i}$ and $V_{j}$ when all the cuts up to time $t_{\varnothing}$ have been removed from $\mathcal T$. On the right, the completion $\mathcal T_{0}$ of $\hat{\cT}_{0}$. Observe that it contains the open path $\rrb y^{i, j}_{00}, y^{i, j}_{01}\llb$, which is isometric to $\rrb V_{i}, x_{\varnothing}\llb$. }
\end{figure}

Let $(i, j)\in \N^{2}$ be  a pair of distinct integers. As explained above, we will introduce a collection $\mathcal R(i, j)$ which will serve as the routing for $(i, j)$. To do so, we need an accompanying collection $\mathcal Y(i, j)=(y^{i, j}_{u})_{u\in \mathbb U^{\ast}}$, in order to record how $\llb V_{i}, V_{j}\rrb $ is fragmented into smaller and smaller subpaths during the pruning process. 
In this part, we explain how to define $y^{i, j}_{u}$ and $r^{i, j}_{u}$ for those $u\in \{0, 1\}^{2}$.  
We will see in Section~\ref{sec: rore} how this construction extends to the remaining $u\in \mathbb U^{\ast}$ in an inductive way. 
From now on, we will fix the pair $(i, j)$ and drop the superscript $^{i, j}$ to ease the notation. 

\smallskip
\noindent
{\bf Generation $0$: } set $y_{0}=V_{i}$ and $y_{1}=V_{j}$; define also $r_{0}=i, r_{1}=j$ and $b_{\varnothing}=r_{0}\wedge r_{1}$. According to the way $\phi$ is defined in Section~\ref{sec: phi-def}, $\phi(b_{\varnothing})\in \cP_{\infty}$ is the location of the first cut on $\llb y_{0}, y_{1}\rrb$. 

\smallskip
\noindent
{\bf Generation $1$: on the side of $\cT$. } In what follows, $q$ stands for a generic binary digit, i.e.~$q=0$ or $q=1$. 
Let $(t_{\varnothing}, x_{\varnothing})$ be the first cut on $\llb y_{0}, y_{1}\rrb$. Let $\hat \cT_{q}$ be the connected component of $\cT\setminus \cP_{t_{\varnothing}}$ that contains $y_{q}$. 
Let $\cT_{q}$ be the completion of $(\hat\cT_{q}, d)$. 
In particular, the closure of $\llb y_{q}, x_{\varnothing}\llb$ in $\cT_{q}$ contains a single element of $\cT_{q}\setminus \hat\cT_{q}$ (i.e.~the ``gap'' left by the removal of $x_{\varnothing}$). We denote by $y_{q1}$ this single element; see Fig.~\ref{fig4} for an example. 
With a slight abuse of notation, we still denote by $\mu$ the sub-probability measure on $\cT_{q}$ defined by $\mu(B)= \mu(B\cap \hat\cT_{q})$ for all Borel sets $B$. 
Set also $y_{q0}=y_{q}$. 
The pruning of $\cT$ extends naturally to a pruning of $\cT_{q}$: the cuts in this case are those $(t_{i}, x_{i})\in\cP$ with $t_{i}\ge t_{\varnothing}$ and $x_{i}\in \hat\cT_{q}$. It then makes sense to define $\cT_{t}(y_{q1})$ to be 
the connected component of $\cT_{q}\setminus \cP_{t}$ containing $y_{q1}$. 
It will be convenient to have this quantity defined for all $t$: we set $\cT_{t}(y_{q1})= \cT_{t}(y_{q})$ if $t<t_{\varnothing}$, 
noting that $x_{\varnothing}$ serves as a natural ``ancestor'' of $y_{q1}$ and it is connected to $y_{q}$ before time $t_{\varnothing}$. 
For $k\ge 1$, let $\tau_{q1, k}$ be the first time $t$ when $\mu(\cT_{t}(y_{q1}))\ne \mu(\cT_{t}(V_{k}))$, which a.s.~corresponds to the moment at which $y_{q1}$ and $V_{k}$ get separated, with the understanding that $y_{q1}=x_{\varnothing}$ for $t<t_{\varnothing}$. 
We now have the proper form of \eqref{dist: delta''}:  for each $k\ge 1$, let
\begin{equation}
\label{def: dist-trace}
\delta_{q1}(k) =  \int_{\tau_{q1, k}}^{\infty} \Big\{\mu\big(\cT_{t}(y_{q1})\big)+\mu\big(\cT_{t}(V_{k})\big)\Big\}dt  \in [0, +\infty], \quad q\in \{0, 1\}.   
\end{equation}

\begin{lem}
\label{prop: trace-def}
Let $q\in \{0, 1\}$. With probability $1$, $(\delta_{q1}(k))_{k\ge 1}$ satisfies (A1-A3) for $(\cC, \delta)$ and the dense sequence $\N$. 
\end{lem}

\noindent
{\bf Generation $1$: on the side of $\cC$. }
Thanks to Lemma~\ref{prop: trace-def} and Lemma~\ref{lem: char}, we see that $(\delta_{q1}(k))_{k\ge 1}$ defines a random point of $\cC$, which we denote as $r_{q1}$,  that verifies $\delta(r_{q1}, k) = \delta_{q1}(k)$ for all $k\in \N$. 
Recall the notation $b_{\varnothing}=r_{0}\wedge r_{1}$. 
Let $\cC^{q}_{b_{\varnothing}}$ be the component of $\cC\setminus\{b_{\varnothing}\}$ that contains $r_{q}$. We also define $r_{q0}=r_{q}$ and $b_{q}=r_{q0}\wedge r_{q1}$. 

\begin{lem}
\label{prop: trace-dist}
Given $\cC^{0}_{b_{\varnothing}}$ and $\cC^{1}_{b_{\varnothing}}$, the collection $\{r_{u}: u\in \{00, 01, 10, 11\}\}$ is independent, and $r_{q1}$ has the same law as $r_{q0}$, $q\in \{0, 1\}$. 
Moreover, up to isometries from $\cT_{q}$ to $\cT$, we have
\[
\rrb y_{0}, y_{1}\llb \;= \;\rrb y_{00}, y_{01}\llb \;\cup\; \rrb y_{10}, y_{11}\llb\; \cup\; \{\phi(b_{\varnothing})\}.
\]
For $q\in \{0, 1\}$, $\phi(b_{q})$ is the location of the first cut on $\rrb y_{q0}, y_{q1}\llb$. 
\end{lem}

Proof of Lemmas~\ref{prop: trace-def} and \ref{prop: trace-dist} will be given in Section~\ref{sec: pf-prop}, after we explain  the routings for $\bp$-trees in Section~\ref{sec: rou-ptree}. 

\subsection{Routings in $\bp$-trees}
\label{sec: rou-ptree}

As we will rely on weak convergence arguments later, it is worthwhile to recall the convergence of $\bp$-trees to ICRT. 
It is common to state the convergence to continuum random trees as measured metric spaces. Therefore we will employ the notions of the Gromov--Hausdorff topology and the Gromov--Prokhorov topology in their standard forms and their pointed versions. We refer to Sections 2.2-2.3 in \cite{BW} for their definitions. 
For each $n\ge 1$, suppose that we have a probability measure $\bp^n=(p^{n}_{1}, p^{n}_{2}, \dots, p^{n}_{n})$ satisfying $p^{n}_{1}\ge p^{n}_{2}\ge \cdots \ge p^{n}_{n}>0$. 
Recall from \eqref{def: pi} the probability distribution $\pi_{\bp^{n}}$ on $\mathbb T_{n}$, the set of all rooted trees labelled by $[n]$. 
Let $T^{n}$ be a random tree of law $\pi_{\bp^{n}}$. 
For $n\ge 1$, denote by $\alpha_{n}$ the positive number such that 
\[
\alpha_{n}^{2}=\sum_{k=1}^{n}(p^{n}_{k})^{2}\,. 
\]
Recall from \eqref{def: Theta} the sequence $\bth=(\theta_1,\theta_2,\dots)$.  Thanks to the scaling property \eqref{id: icrt-scale}, we can always assume $\|\bth\|\le 1$. Let $\beta=1-\|\bth\|^{2}$. 
Suppose that 
\begin{equation}\label{p-tree}
\lim_{n\to\infty}\alpha_n=0,
\qquad\text{and} \qquad
\lim_{n\to\infty}\frac{p^{n}_{k}}{\alpha_n}=\theta_i, \quad \text{ for each }k\ge 1.
\end{equation}
Denote by  $(T^{n}, d_n)$ the metric space obtained by endowing the vertex set $[n]$ with 
the graph distance $d_n$ of $T^{n}$. 
Then, with $(\cT,d,\mu)$ an ICRT under $\bP^{\beta, \bth}$, Aldous \& Pitman \cite{Al00} show that
\begin{equation}\label{eq: CP}
(T^n, \,\alpha_{n}\,d_n,\, \bp^n) \xrightarrow[n\to\infty]{(d)} (\cT,\, d, \,\mu),
\end{equation}
with respect to the Gromov--Prokhorov topology. Moreover, as shown in \cite{BW}, the pruning 
process on $\cT$ described in Section~\ref{sec:pruning_cut-trees} happens to be the scaling limit 
of the following discrete pruning process on $T^{n}$. 

\medskip
\noindent\textbf{Pruning $\bp$-trees and the discrete cut tree.}
Let $\cP^{n}=\{(t^{n}_{i}, x^{n}_{i}): 0\le t^{n}_{1}< t^{n}_{2}< \cdots\}$ be a Poisson point 
process on $\R_{+}\times  T^{n}$ with intensity $dt \otimes\bp^{n}(dx)$.  Just as in the 
continuous version, we remove the vertex $x^{n}_{i}$ (if it has not yet been removed) at time 
$t^{n}_{i}$. Note that distribution-wise this amounts to the same pruning introduced in Section~\ref{sec: intro-p-tree}. 
Recall that we have already defined in Section~\ref{sec: intro-p-tree} a discrete cut tree $C^{n}$ for this pruning of $T^{n}$. 

Let $\delta^n$ denote the graph distance on $C^n$, and  suppose that \eqref{p-tree} holds. 
Recall that $C^{n}$ is rooted at the first removed vertex in $T^{n}$; let us denote it by $\hat\rho^{n}$. 
Recall also that $0$ is the root of $\cC$. 
The following convergence shows that the cut tree $\cC$ can be seen as the continuous counterpart 
of discrete tree $C^{n}$: we have (Theorem~3.4 in \cite{BW})
\begin{equation}
\label{cv: cut-tree}
(C^n, \,\alpha_{n}\, \delta^n,\, \bp^n,\, \hat\rho^{n}) \xrightarrow[n\to\infty]{(d)} (\cC,\, \delta, \,\nu,\, 0),
\end{equation}
with respect to the pointed Gromov--Prokhorov topology, and jointly with the convergence in \eqref{eq: CP}. 

Rather remarkably, $C^{n}$ is also distributed as a $\bp_{n}$-tree (Theorem 4.10 in \cite{BW}). 
This explains the analogue equality in distribution between $\cT$ and $\cC$ in the limit: 
combined with \eqref{cv: cut-tree}, this implies that 
$(\cC, \delta, \nu)$ is a $(\beta, \bth)$-ICRT under $\bP^{\beta, \bth}$. 

\medskip
\noindent\textbf{Routings in $T^{n}$ and $C^{n}$. } 
Let $(V^{n}_{i})_{i\ge 0}$ be a sequence of i.i.d~points of common distribution $\bp^{n}$, independent of $\cP^{n}$. 
Set $y^{n}_{0}=r^{n}_{0}=V^{n}_{i}$ and $y^{n}_{1}=r^{n}_{1}=V^{n}_{j}$; let $(t^{n}_{\varnothing}, b^{n}_{\varnothing})$ be the first cut on $\llb y^{n}_{0}, y^{n}_{1}\rrb$. 

We use the notation $m|_{B}$ to denote the restriction of a finite measure $m$ to a subset $A$, namely, $m|_{A}(B)=m(A\cap B)/m(A)$ for all measurable sets $B$, as long as $m(A)>0$. 

\begin{lem}
\label{lem: p-tree}
Let $q\in \{0, 1\}$. Let $B_{0}, B_{1}$ be two disjoint non-empty subsets of $[n]$. 
Let $T^{n}_{q}$ be the subtree of $T^{n}$ on the vertex set 
\[
\mathcal V(T^{n}_{q})=\{v\in T^{n}: \mathcal P^{n}\cap [0, t^{n}_{\varnothing}]\times \llb y^{n}_{q}, v\rrb=\varnothing\},
\]
namely the subtree formed by the vertices connected to $y^{n}_{q}$ right after the cut at $b^{n}_{\varnothing}$. If $\mathcal V(T^n_q)$ is non-empty, which occurs if $b^{n}_{\varnothing}\ne y^{n}_{q}$, we let $T^n_q$ be rooted at $y^{n}_{q}$. 
The following statements hold true.
\begin{enumerate}[(i)]
	\item
	$b^{n}_{\varnothing}$ is the most recent common ancestor of $r_1^{n}=y^{n}_{1}$ and  $r_0^{n}=y^{n}_{0}$ in $C^{n}$. 
	\item
	Let $C^{n, q}_{b^{n}_{\varnothing}}$ be the connected subgraph of $C^n\setminus\{b^{n}_{\varnothing}\}$ that contains $r^n_q$. Then $C^{n, q}_{b^{n}_{\varnothing}}$ is the cut tree of $T^{n}_{q}$ with respect to the cuts from $\cP^{n}$ restricted to $T^{n}_{q}$; in particular $C^{n, q}_{b^{n}_{\varnothing}}$ and $T^n_q$ have the same vertex set. 
	
	\item Let $y^n_{q1}$ be the unique vertex of $C^{n, q}_{b^{n}_{\varnothing}}$ that was adjacent to $b^n_{\varnothing}$ in $T^n$. 
	Conditionally on the vertex sets of $C^{n, 0}_{b^{n}_{\varnothing}}$ and $C^{n, 1}_{b^{n}_{\varnothing}}$ being $B_0$ and $B_1$ respectively,  
	\begin{itemize}
		\item $(y^n_0, y^n_{01}, y^{n}_{1}, y^{n}_{11})$ is independent;
		\item $y^{n}_{q}$, $y^{n}_{q1}$ both have distribution $\bp^n|_{B_q}$ in $C^{n, q}_{b^{n}_{\varnothing}}$.
	\end{itemize}
	The same properties also hold when replacing $C^{n, q}_{b^{n}_{\varnothing}}$ by $T^{n}_{q}$. 
\end{enumerate}
\end{lem}

The above says that $T^{n}_{0}$ and $T^{n}_{1}$ each contain a neighbour  of $b^{n}_{\varnothing}$, namely $y^{n}_{01}$ and $y^{n}_{11}$, whenever $T^{n}_{0}$ and $T^{n}_{1}$ are non empty. Their copies in $C^{n}$  serve as the discrete counterparts of $r_{01}, r_{11}$. We thus set $r^{n}_{q1}=y^{n}_{q1}$, $q= 0, 1$.   

\begin{proof}[Proof of Lemma~\ref{lem: p-tree}] 
We rely on the results of \cite{BW}, where this is already implicit. 
Note that \emph{(i)} and \emph{(ii)} are straightforward from the definition of the cut tree in there. 
It remains to see \emph{iii)}. 
By Theorem~4.10 of \cite{BW}, the cut tree $C^n$ is a $\bp^n$-tree, and $V^n_i$, $V^n_j$ are two independent nodes with distribution $\bp^n$. Therefore, given the vertex set being $B_0$ (resp.~$B_1$), 
it can be checked from the product form of $\pi_{\bp^{n}}$ that $y^{n}_{0}=V^{n}_{i}$ (resp.~$y^{n}_{1}=V^{n}_{j}$) has distribution $\bp^{n}|_{B_{0}}$ inside $C^{n, 0}_{b^{n}_{\varnothing}}$ (resp.~$\bp^{n}|_{B_{1}}$ inside $C^{n, 1}_{b^{n}_{\varnothing}}$). Since $T^{n}_{q}$ shares the same vertex set as $C^{n, q}_{b^{n}_{\varnothing}}$, the same can be said 
for $y^{n}_{q}$ inside $T^{n}_{q}$. 
Consider now the discrete reconstruction procedure described on p.2406-2407 in \cite{BW}. In particular, the construction there implies that $y^{n}_{q1}$ is sampled using $\bp^{n}$ conditioned to be in $C^{n, q}_{b^{n}_{\varnothing}}$. Moreover, when conditioned on the vertex sets, these samplings are  independent of each other, and of $y^{n}_{0}, y^{n}_{1}$. This proves the statement for $C^{n}$. Since $T^{n}_{q}$ shares the same vertex set as $C^{n, q}_{b^{n}_{\varnothing}}$, the distributions of $y^{n}_{q}, y^{n}_{q1}$ are unchanged. 
\end{proof}

To prepare for our proof of Lemmas~\ref{prop: trace-def} and \ref{prop: trace-dist}, let us note some useful properties of $C^{n}$. 
We have denoted by $\delta^{n}$ the graph distance in $C^{n}$. Recall also that $\hat\rho^{n}$ stands for the root of $C^{n}$. 
For a vertex $v$ of $T^{n}$ and $t\ge 0$, we write $T^{n}_{t}(v)$ for the subgraph of $T^{n}$ formed by vertices still connected to $v$ at time $t$. 
As a consequence of our definition of $C^{n}$, the distance $\delta^{n}$ has the following interpretation in terms of  the pruning of $T^{n}$: for all $v\in [n]$, we have
\begin{equation}
\label{eq: record}
\delta^{n}(\hat\rho^{n}, v) + 1 =\#\big\{(t^{n}_{i}, x^{n}_{i})\in \cP^{n}: T^{n}_{t^{n}_{i}}(v)\subsetneq T^{n}_{t^{n}_{i}-}(v)\big\} \overset{\text{a.s.}}{=} \#\big\{(t^{n}_{i}, x^{n}_{i})\in \cP^{n}: x^{n}_{i}\in T^{n}_{t^{n}_{i}-}(v)\big\}.
\end{equation}
In words, $\delta^{n}(\hat\rho^{n}, v) + 1$ counts the number of jumps of $(T^{n}_{t}(v))_{t\ge 0}$.  
Alternatively, this can also be formulated as follows: say a vertex $u$ is a {\it record for $v$} if it is the first on its path to $v$ to be removed during the pruning of $T^{n}$; then $\delta^{n}(\hat\rho^{n}, v) + 1$ corresponds to the number of records for $v$. Since $\delta^{n}$ satisfies that for any $u, v\in [n]$, 
\[
\delta^{n}(u, v)=\delta^{n}(\hat\rho^{n}, u)+\delta^{n}(\hat\rho^{n}, v)-2\delta^{n}(\hat\rho^{n}, u\wedge v),
\]
where $u\wedge v$ is the most recent common ancestor of $u$ and $v$ in $C^{n}$, we can extend \eqref{eq: record} to any pair of vertices $u, v\in [n]$:
\begin{align}\notag
\delta^{n}(u, v) = &\, \#\big\{(t^{n}_{i}, x^{n}_{i})\in \cP^{n}: t^{n}_{i}> \tau^{n}(u, v), T^{n}_{t^{n}_{i}}(u)\subsetneq T^{n}_{t^{n}_{i}-}(u)\big\} \\ \label{eq: dist-record}
& + \#\big\{(t^{n}_{i}, x^{n}_{i})\in \cP^{n}: t^{n}_{i}> \tau^{n}(u, v), T^{n}_{t^{n}_{i}}(v)\subsetneq T^{n}_{t^{n}_{i}-}(v)\big\},
\end{align}
where $\tau^{n}(u, v)=\inf\{t: T^{n}_{t}(u)\ne T^{n}_{t}(v)\}$ is the moment when $u$ is separated from $v$. 

\subsection{Proof of Lemmas~\ref{prop: trace-def} and \ref{prop: trace-dist} via weak convergence arguments}
\label{sec: pf-prop}

Recall the quantities $\delta_{u}(k), k\in \N, u\in \{01, 11\}$ defined in \eqref{def: dist-trace}. 
Recall the vertices $r^{n}_{01}=y^{n}_{01}, r^{n}_{11}=y^{n}_{11}$ introduced in Lemma~\ref{lem: p-tree}, as well as  the graph distance $\delta^{n}$ of the discrete cut tree $C^{n}$. 
Lemmas~\ref{prop: trace-def} and \ref{prop: trace-dist} are consequences of the following
\begin{lem}
\label{prop: cvdelta}
We have $\delta_{u}(k)<\infty$ a.s.~for each $k\in \N$ and $u\in \{01, 11\}$. 
Let $(V^{n}_{k})_{k\ge 0}$ be a sequence of i.i.d.~vertices of common law $\bp^{n}$. 
Under the assumption in \eqref{p-tree}, 
jointly with the convergences in \eqref{eq: CP} and \eqref{cv: cut-tree}, we have 
\begin{equation}
	\label{cv: Cdist}
	\Big(\alpha_{n}\delta^{n}(V^{n}_{k}, r^{n}_{u}): k\in \N, u\in\{01, 11\}\Big) \xrightarrow[n\to\infty]{(d)} \Big(\delta_{u}(k): k\in \N, u\in\{01, 11\}\Big),
\end{equation}
with respect to the product topology of $\R^{\N}\times\R^{\N}$. Moreover, given $\cC^{q}_{b_{\varnothing}}$, $(\delta_{q1}(k))_{k\ge 1}$ has the same law as $(\delta({i}, U'))_{i\ge 1}$, where $U'$ is an independent point with distribution $\nu|_{\cC^{q}_{b_{\varnothing}}}$, $q\in\{0, 1\}$. 
\end{lem}

\begin{proof}
This follows closely the arguments in Section 5.2 of \cite{BW}; therefore we only outline the main steps. Recall the spanning tree $\cR'_{k}=\cup_{1\le j\le k}\llb V_{0}, V_{j}\rrb$ of $\cT$. Define similarly its discrete counterparts: let $R^{n}_{k}=\cup_{1\le j\le k}\llb V^{n}_{0}, V^{n}_{j}\rrb$ be the subtree of $T^{n}$ spanned by $V^{n}_{0}, \dots, V^{n}_{k}$. Then Eq.~(5.16) in \cite{BW} says that for all $k\ge 1$ and $t\ge 0$,
\begin{equation}
	\label{pfeq: cv-cut}
	\big(R^{n}_{k}, \alpha_{n}d_{n}, \cP^{n}\cap [0, t/\alpha_{n}]\times R^{n}_{k}\big) \xrightarrow[n\to\infty]{(d)} \big(\cR'_{k}, d, \cP\cap [0, t]\times \cR'_{k}\big)
\end{equation}
with respect to the pointed Gromov--Hausdorff topology. Let $q\in \{0, 1\}$. Note that if $V_{m}\in \cT_{q}$, then $d_{\cT_{q}}(V_{m}, y_{q1}) = d_{\cT}(V_{m}, x_{\varnothing})$, where we recall that $x_{\varnothing}$ is the location of the first cut on $\llb y_{0}, y_{1}\rrb$. We deduce that

\noindent
{\bf Step 1.} For each $k\ge 1$ and $q\in \{0, 1\}$, we have
\begin{equation}
	\label{pfcv: trace}
	\big(R^{n}_{k}\cap T^{n}_{q}, \alpha_{n}d_{n}, y^{n}_{q1}\big) \xrightarrow[n\to\infty]{(d)} \big(\cR'_{k}\cap \cT_{q}, d, y_{q1}\big)
\end{equation}
with respect to the pointed Gromov--Hausdorff topology.

\medskip
\noindent
{\bf Step 2.} By adapting the arguments in the proof of Lemma 5.4 in \cite{BW}, we can show that 
\begin{equation}
	\label{pfcv: mass}
	\Big(\bp^{n}\big(T^{n}_{t/\alpha_{n}}(y^{n}_{q1})\big): t\ge 0, q\in \{0, 1\}\Big) \xrightarrow[n\to\infty]{(d)} \Big(\mu\big(\cT_{t}(y_{q1})\big): t\ge 0, q\in \{0, 1\}\Big)
\end{equation}
uniformly on all compact sets of $\R_{+}$, and jointly with 
\begin{equation}
	\label{pfcv: sep-time}
	\big(\alpha_{n}\tau^{n}(V^{n}_{k}, y^{n}_{q1}) : k\ge 1, q\in \{0, 1\}\big) \xrightarrow[n\to\infty]{(d)} \big(\tau_{q1, k}: k\ge 1, q\in \{0, 1\}\big).
\end{equation}


\medskip
\noindent
{\bf Step 3.} Recall the notion of record from the paragraph below \eqref{eq: record}. For any vertex $v$, we introduce $\tilde\delta^{n}_{t}(v)$ as the number of records for $v$ that occur before time $t$. Following the proof of Lemma 5.6 in \cite{BW}, we can show that 
\[
\tilde\delta^{n}_{t}(y^{n}_{q})-\int_{0}^{t}\bp^{n}\big(T^{n}_{s}(y^{n}_{q1})\big)ds, \quad t\ge 0,
\]
is a martingale. Moreover, we have 
\begin{align}
	\notag
	\mathbb E\Big[\sup_{t\ge 0}\Big(\tilde\delta^{n}_{t}(y^{n}_{q})-\int_{0}^{t}\bp^{n}\big(T^{n}_{s}(y^{n}_{q1})\big)ds \Big)^{2}\Big]&\le 4\mathbb E\Big[\int_{0}^{\infty}\bp^{n}\big(T^{n}_{s}(y^{n}_{q1})\big)\Big] \\ \label{pf: mart}
	&= 4 \mathbb E[\tilde\delta^{n}_{\infty}(y^{n}_{q})] = 4\mathbb E[\delta^{n}(\hat\rho^{n}, y^{n}_{q})+1],
\end{align}
where we have used the martingale property for the first identity and \eqref{eq: record} for the second. 
Let $U^{n}$ be a random vertex with law $\bp^{n}$. Then according to Lemma~\ref{lem: p-tree}, $y^{n}_{q}$ is distributed as $U^{n}$ conditioned to be in $T^{n}_{q}$. In particular, 
\[
\bp^{n}(T^{n}_{q})\cdot \mathbb E[\delta^{n}(\hat\rho^{n}, y^{n}_{q})\,|\,T^{n}] = \mathbb E[\delta^{n}(\hat\rho^{n}, U^{n})\,|\,T^{n}].
\]
Note that $\sup_{n\ge 1}\alpha_{n}\mathbb E[\delta^{n}(\hat\rho^{n}, U^{n})]<\infty$ (see for instance Eq.~(4) in \cite{Al00}), while \eqref{pfcv: trace} implies that $\bp^{n}(T^{n}_{q})\overset{(d)}{\to} \mu(\cT_{q})>0$, we infer that 
\[
\sup_{n\ge 1}\mathbb E\Big[\int_{0}^{\infty}\bp^{n}\big(T^{n}_{s/\alpha_{n}}(y^{n}_{q1})\big)ds\Big]=\sup_{n\ge 1}\alpha_{n}\mathbb E\Big[\int_{0}^{\infty}\bp^{n}\big(T^{n}_{s}(y^{n}_{q1})\big)ds\Big]<\infty.
\]
Combined with \eqref{pfcv: mass}, applying monotone convergence and  Fatou's Lemma leads to
\begin{equation}
	\label{pf: finite}
	\mathbb E\Big[\int_{0}^{\infty}\mu\big(\cT_{s}(y_{q1})\big)ds\Big]=\lim_{M\to\infty}\mathbb E\Big[\int_{0}^{M}\mu\big(\cT_{s}(y_{q1})\big)ds\Big]\le \liminf_{n\to\infty}\mathbb E\Big[\int_{0}^{\infty}\bp^{n}\big(T^{n}_{s/\alpha_{n}}(y^{n}_{q1})\big)ds\Big]<\infty.
\end{equation}

\medskip
\noindent
{\bf Step 4.} With similar arguments as before and Lemma~5.8 in \cite{BW}, we can show that for all $\epsilon>0$, \begin{equation}
	\label{pf: tight}
	\lim_{t\to\infty}\limsup_{n\to\infty}\mathbb P\Big(\delta^{n}(\hat\rho^{n}, y^{n}_{q1})-\tilde\delta^{n}_{t/\alpha_{n}}(y^{n}_{q1})>\epsilon/\alpha_{n}\Big) = 0.
\end{equation}

\medskip
\noindent
{\bf Step 5.} Combining \eqref{pfcv: mass}, \eqref{pf: mart}, \eqref{pf: finite} and \eqref{pf: tight}, we deduce  the following joint convergence:
\begin{align*}
	\big(\alpha_{n}\delta^{n}(\hat\rho^{n}, y^{n}_{q1})\big)_{q\in \{0,1\}} &\xrightarrow[n\to\infty]{(d)} \Big(\int_{0}^{\infty}\mu(\cT(y_{q1}))dt\Big)_{q\in \{0, 1\}} \\
	\big(\alpha_{n}\tilde\delta^{n}_{t/\alpha_{n}}(y^{n}_{q1}): t\ge 0\big)_{q\in \{0,1\}} &\xrightarrow[n\to\infty]{(d)} \Big(\int_{0}^{t}\mu(\cT(y_{q1}))dt: t\ge 0\Big)_{q\in \{0, 1\}} \quad \text{uniformly in $t$},
\end{align*}
jointly with the convergence in \eqref{pfcv: sep-time}. Together with \eqref{eq: dist-record}, this implies \eqref{cv: Cdist}. Finally, by Lemma~\ref{lem: p-tree}, $y^{n}_{q1}\sim \bp^{n}|_{C^{n, q}_{b^{n}_{\varnothing}}}$ in $C^{n}$. It follows that the law of the sequence
\[
\Big((C^{n}, \delta^{n}, \bp^{n}, r^{n}_{0}, r^{n}_{1}),\, (C^{n, q}_{b^{n}_{\varnothing}}, \delta^{n}, \bp^{n}|_{C^{n, q}_{b^{n}_{\varnothing}}}),\, (\delta^{n}(V^{n}_{k}, y^{n}_{q1}))_{k\in \N}\Big), \, n\in \N, 
\] 
has a unique limit point, which is the law of $((\cC, \delta, \nu), (\cC^{q}_{b_{\varnothing}}, \delta, \nu|_{\cC^{q}_{b_{\varnothing}}}), (\delta({k}, U'))_{k\in \N})$, with $U'\sim \nu|_{\cC^{q}_{b_{\varnothing}}}$ and independent of $r_{q}$. This completes the proof. 
\end{proof}

\begin{proof}[Proof of Lemma~\ref{prop: trace-def}]

By Lemma~\ref{prop: cvdelta}, $(\delta_{q1}(i))_{i\ge 1}$ has the same law as $(\delta(i, U'))_{i\ge 1}$, where $U'$ is some random point in $\cC$. It follows that Properties (A1-A3) are satisfied by $(\delta_{q1}(i))_{i\ge 1}$ a.s.
\end{proof}

We will need the following result for the proof of Lemma~\ref{prop: trace-dist}. 

\begin{lem}
Let $q\in \{0, 1\}$. 
We have a.s.
\begin{equation}
	\label{id: dist-to-root}
	\delta(r_{q1}, 0) = \int_{0}^{\infty}\mu\big(\cT_{t}(y_{q1})\big)dt = \int_{0}^{t_{\varnothing}}\mu\big(\cT_{t}(y_{q})\big)dt + \int_{t_{\varnothing}}^{\infty}\mu\big(\cT_{t}(y_{q1})\big)dt. 
\end{equation}
\end{lem}

\begin{proof}
We will only give detailed arguments for $q=0$, as the other case is similar. 
Since $\N$ is dense in $\cC$, we can find a sequence of random points $(I_{n})_{n\in \N}$ which converges to $r_{01}$ a.s. Recall that $\tau_{01, i}$ is the separation time of $y_{01}$ from $V_{i}$ and that we have $\cT_{t}(y_{01})=\cT_{t}(V_{i})$ for all $t<\tau_{01, i}$.  It follows that
\begin{multline} 
	\left|\int_{0}^{\infty}\mu\big(\cT_{t}(y_{01})\big)dt - \int_{0}^{\infty}\mu\big(\cT_{t}(V_{I_{n}})\big)dt\right| \\
	=  \left|\int_{\tau_{01, I_{n}}}^{\infty}\mu\big(\cT_{t}(y_{01})\big)dt - \int_{\tau_{01, I_{n}}}^{\infty}\mu\big(\cT_{t}(V_{I_{n}})\big)dt\right| \le \delta(r_{01}, I_{n})\to 0.
\end{multline}
As $\delta(0, r_{01})=\lim_{n\to\infty}\delta(0, I_{n})$, the conclusion follows. 
\end{proof}

\begin{proof}[Proof of Lemma~\ref{prop: trace-dist}]
The first statement is a consequence of the last statement in Lemma~\ref{prop: cvdelta}. 
The second one is immediate as $\phi(b_{\varnothing})$ is the location of the first cut on $\llb y_{0}, y_{1}\rrb$, and from the definitions of $(y_{u}), u\in \{0, 1\}^{2}$. 
For the third statement, we note that, by the definition of $\cC$, there is some $(t_{\ast}, x_{\ast})\in \cP$ so that $x_{\ast}=\phi(b_{q})$. 
Let $\cC^{0}_{b_{q}}$ (resp.~$\cC^{1}_{b_{q}}$) be the subtree of $\cC$ above $b_{q}$ that contains $r_{q0}$ (resp.~$r_{q1}$). Since $\N$ is dense in $\cC$, we can find some $j_{0}\in \cC^{0}_{b_{q}}$ and $j_{1}\in \cC^{1}_{b_{q}}$. Note that we have $b_{q}=j_{0}\wedge j_{1}=r_{q0}\wedge j_{1}=j_{0}\wedge r_{q1}$. Then by the construction of $\cC$, $(t_{\ast}, x_{\ast})$ is the first cut on $\llb V_{j_{0}}, V_{j_{1}}\rrb$. Let us show it is also the first cut on $\llb y_{q0}, y_{q1}\rrb$. To that end, let $b'=j_{1}\wedge r_{q1}$, which belongs to $\cC^{1}_{b_{q}}$. It follows that $\delta(r_{q1}, b')<\delta(r_{q1}, b_{q})$. Recall that $\tau_{q1, k}$ stands for the time when $y_{q1}$ get separated from $V_{k}$, with the convention that $y_{q1}=x_{\varnothing}$ if $t< t_{\varnothing}$. We claim that 
\[
\tau_{q1, j_{1}} > \tau_{q1, j_{0}}.
\]
To see why this is true, recall that $\delta(r_{q1}, j_{0}) = \delta_{q1}(j_{0})$ is given in \eqref{def: dist-trace}. Since $\delta$ is the distance function of an $\R$-tree, it verifies 
\[
2\delta(r_{q1}, b_{q}) = 2\delta(r_{q1}, r_{q1}\wedge j_{0}) = \delta(r_{q1}, j_{0})+\delta(r_{q1}, 0)-\delta(j_{0}, 0). 
\]
Plugging \eqref{id: dist-to-root} into the right-hand side, we find that 
\[
\delta(r_{q1}, b_{q}) = \int_{\tau_{q1, j_{0}}}^{\infty}\mu\big(\cT_{t}(y_{q1})\big)dt.
\] 
A similar argument yields
\[
\delta(r_{q1}, b') = \int_{\tau_{q1, j_{1}}}^{\infty}\mu\big(\cT_{t}(y_{q1})\big)dt.
\]
Comparing this with $\delta(r_{q1}, b')<\delta(r_{q1}, b_{q})$, we deduce that $\tau_{q1, j_{1}} > \tau_{q1, j_{0}}$. In consequence, we have $y_{q1}\in \cT_{t}(V_{j_{1}})$ for all $t\le t_{\ast}$. Replacing $y_{q1}$ with $y_{q0}$ (and \eqref{def: delta}, \eqref{def: delta'} substituting for \eqref{def: dist-trace}, \eqref{id: dist-to-root}), we find that $y_{q0}\in \cT_{t}(V_{j_{0}})$ for all $t\le t_{\ast}$. As $(t_{\ast}, x_{\ast})$ is the first cut on $\llb V_{j_{0}}, V_{j_{1}}\rrb$, it follows $\cT_{t}(V_{j_{0}})=\cT_{t}(V_{j_{1}})$ for $t<t_{\ast}$. 
This means that $y_{q0}$ is connected to $y_{q1}$ right up to time $t_{\ast}$ but becomes disconnected at the moment $t_{\ast}$. Therefore, $(t_{\ast}, x_{\ast})$ is also the first cut on $\llb y_{q0}, y_{q1}\rrb$. The conclusion follows. 
\end{proof}

\subsection{Constructing $\mathcal R(i, j)$: further generations}
\label{sec: rore}

Having built the first generation of the collection $\mathcal R(i, j)$, we simply iterate the previous procedure to define the rest. More precisely, suppose that we have the collections:
\[
\big\{y_{u}: u\in \{0, 1\}^{m}, 1\le m\le n\big\}, \, \big\{r_{u}: u\in \{0, 1\}^{m}, 1\le m\le n\big\}, \,
\big\{\cT_{u}: u\in \{0, 1\}^{m}, 1\le m\le n\big\},
\]
as well as 
\[
\{b_{u}: u\in \{0, 1\}^{m-1}, 1\le m\le n\}.
\]
In above, each $r_{u}$ is a point of $\cC$ and $b_{u}$ is the most recent common ancestor of $r_{u0}$ and $r_{u1}$ in $\cC$. Moreover, $\cT_{u}$ is a complete real tree containing a subset $\hat{\cT}_{u}\subseteq \cT$ which satisfies that  $\cT_{u}\setminus \hat{\cT}_{u}$ is at most countable and $y_{u0}, y_{u1}\in \cT_{u}$. As a consequence, we can introduce a pruning of $\cT_{u}$ by using the cuts $(t_{i}, x_{i})\in \cP$ satisfying $x_{i}\in \hat\cT_{u}$. 

\smallskip
\noindent
{\bf Inductive definition for the routings. }
Let $u\in\{0, 1\}^{n-1}$ and $q\in \{0, 1\}$. 
Let $(t_{u}, x_{u})$ be the first cut on $\llb y_{u0}, y_{u1}\rrb$ in the pruning of $\cT_{u}$. Let $\hat \cT_{uq}$ be the connected component of $\cT_{u}\setminus \cP_{t_{u}}$ that contains $y_{uq}$. Let $\cT_{uq}$ be the completion of $(\hat\cT_{uq}, d)$ and set $y_{uq1}$ to be the single element in the closure of $\llb y_{uq}, x_{u}\llb$  that does not belong to $\hat\cT_{uq}$. 
Set also $y_{uq0}=y_{uq}$. 
The pruning of $\cT$ extends naturally to a pruning of $\cT_{uq}$: the cuts in this case are those $(t_{i}, x_{i})\in\cP$ with $t_{i}\ge t_{u}$ and $x_{i}\in \hat\cT_{uq}$. We then define $\cT_{t}(y_{uq1})$ to be 
the connected component of $\cT_{uq}\setminus \cP_{t}$ containing $y_{uq1}$, with the convention that  
$\cT_{t}(y_{uq1})= \cT_{t}(y_{uq})$ if $t<t_{u}$. 
We also extend the measure $\mu$ to $\cT_{uq}$ by setting $\mu(B)=\mu(B\cap \hat{\cT_{uq}})$ for all Borel measure $B\subseteq \cT_{uq}$. 
For $k\ge 1$, let $\tau_{uq1, k}$ be the first time $t$ when $\mu(\cT_{t}(y_{uq1}))\ne \mu(\cT_{t}(V_{k}))$. 
We then define for each $k\ge 1$: 
\begin{equation}
\label{id: dist-ro}
\delta_{uq1}(k) =  \int_{\tau_{uq1, k}}^{\infty} \Big\{\mu\big(\cT_{t}(y_{uq1})\big)+\mu\big(\cT_{t}(V_{k})\big)\Big\}dt .  
\end{equation}
Following the same arguments as in the proof of Lemma~\ref{prop: trace-def}, we can show the following result.
\begin{lem}
Let $u\in\{0, 1\}^{n-1}$ and $q\in \{0, 1\}$. With probability $1$, $(\delta_{uq1}(k))_{k\ge 1}$ satisfies (A1-A3) for $(\cC, \delta)$ and the dense sequence $\N$. 
\end{lem}
It follows that $(\delta_{uq1}(k))_{k\ge 1}$ defines a random point $r_{uq1}$ of $\cC$ characterised by the identities $\delta(r_{uq1}, k) = \delta_{uq1}(k)$ for all $k \in \N$. Let $\cC^{q}_{b_{u}}$ be the component of $\cC\setminus\{b_{u}\}$ that contains $r_{uq}$. Set also $r_{uq0}=r_{uq}$ and $b_{uq}=r_{uq0}\wedge r_{uq1}$. 
This completes the definition for $y_{w}, r_{w}$ for all $w\in \{0, 1\}^{n+1}$ and $b_{w'}$ for all $w'\in \{0, 1\}^{n}$.  
Moreover, by adapting the proof of Lemma~\ref{prop: trace-dist}, we can show the following statements also hold true. 
\begin{lem}
\label{prop: trace-dist'}
Let $u\in\{0, 1\}^{n-1}$ and $q\in \{0, 1\}$. Given $\cC^{0}_{b_u}$ and $\cC^{1}_{b_u}$, the collection $\{r_{u00}, r_{u01}, r_{u10}, r_{u11}\}$ is independent, and $r_{uq1}$, $r_{uq1}$ both have the law $\nu|_{\cC^{q}_{b_{u}}}$.  
Moreover, up to isometries from $\cT_{uq}$ to $\cT$, we have
\[
\rrb y_{u0}, y_{u1}\llb \;= \;\rrb y_{u00}, y_{u01}\llb \;\cup\; \rrb y_{u10}, y_{u11}\llb\; \cup\; \{\phi(b_{u})\}.
\]
For $q\in \{0, 1\}$, $\phi(b_{uq})$ is the location of the first cut on $\rrb y_{uq0}, y_{uq1}\llb$. 
\end{lem}

The above allows us to define the collections $\mathcal R(i, j)=\{r^{i, j}_{u}: u\in \mathbb U^{\ast}\}$ and $\mathcal B(i, j)=\{b^{i, j}_{u}: u\in\mathbb U\}$, for all distinct pairs $(i, j)\in \N^{2}$. 

\begin{lem}
\label{lem: ro-br}
We have $b^{i, j}_{u}\in \mathcal B(i, j)$ if and only if $\phi(b^{i, j}_{u})\in \llb V_{i}, V_{j}\rrb$.
\end{lem}

\begin{proof}
Lemma \ref{prop: trace-dist'} ensures that for all $u\in \{0, 1\}^{n}$ and $n\ge 0$, $\rrb y^{i, j}_{u0}, y^{i, j}_{u1}\llb$ is a subset of $\llb V_{i}, V_{j}\rrb$. It also says that $\phi(b^{i, j}_{u})$ is the location of a cut on $\rrb y^{i, j}_{u0}, y^{i, j}_{u1}\llb$. It follows that $\phi(b^{i, j}_{u})\in  \llb V_{i}, V_{j}\rrb$. 
It remains to show the converse is also true. 
We note that the previous construction implies the cut at $\phi(b^{i, j}_{u0})$ occurs after the cut at $\phi(b^{i, j}_{u})$. Together with Lemma~\ref{prop: trace-dist'}, this implies that $\{\phi(b^{i, j}_{u}): u\in \{0, 1\}^{m}, 0\le m<n\}$ contains at least the locations of the first $n$ cuts on $\llb V_{i}, V_{j}\rrb$. Since this is true for all $n$, the proof is complete. 
\end{proof}

\subsection{Assembling the collection of traces}

The collections of routings $\mathcal R(i, j)=\{r^{i, j}_{u}: u\in \mathbb U^{\ast}\}$, $i\ne j\in \N$, have the following consistency property:

\begin{lem}
\label{lem: cons}
Let $(i, j)\in\N^{2}, (k, \ell)\in \N^{2}$ be such that $i\ne j, k\ne \ell$. Suppose that $b^{i, j}_{u}=b^{k, \ell}_{w}$ for some $u, w\in \mathbb U^{\ast}$, and that $r^{i, j}_{u0}$ and $r^{k, \ell}_{w0}$ are in the same component of $\cC\setminus\{b^{i, j}_{u}\}$. Then we have $r^{i, j}_{u1}=r^{k, \ell}_{w1}$. 
\end{lem}

\begin{proof}
Since $r^{i, j}_{u1}$ and $r^{k, \ell}_{w1}$ are characterised by \eqref{id: dist-ro}, it suffices to show that $\mu(\cT_{t}(y^{i, j}_{u1}))=\mu(\cT_{t}(y^{k, \ell}_{w1}))$ for all $t\ge 0$. By our previous construction, there is some $(t_{\ast}, x_{\ast})\in \cP$ with $x_{\ast}=\phi(b^{i, j}_{u})$ that is the first cut on the path from $y^{i, j}_{u}$ to another point $y^{i, j}_{u'}$ with $u'$ in the same generation of $\mathbb U$ as $u$. As $b^{i, j}_{u}=b^{k, \ell}_{w}$, the same pair $(t_{\ast}, x_{\ast})$ is also the first cut on the path from $y^{k, \ell}_{w}$ to $y^{k, \ell}_{w'}$ for some $w'\in \mathbb U$. Following the previous notation, $\hat\cT^{i, j}_{u}$ (resp.~$\hat\cT^{k, \ell}_{w}$) is the connected component  of $\cT\setminus\cP_{t_{\ast}}$ that contains $y^{i, j}_{u}$. By assumption, $r^{i, j}_{u0}$ and $r^{k, \ell}_{w0}$ are in the same component of $\cC\setminus\{b^{i, j}_{u}\}$. Adapting the arguments in the proof of Lemma~\ref{prop: trace-dist}, we can show that $y^{i, j}_{u}$ is separated from $y^{k, \ell}_{w}$ after time $t_{\ast}$. It follows that $\hat\cT^{i, j}_{u}=\hat\cT^{k, \ell}_{u}$, and therefore $\cT^{i, j}_{u}$ is isometric to $\cT^{k, \ell}_{u}$. In particular, this isometry maps $y^{i, j}_{u1}$ to $y^{k, \ell}_{w1}$. It follows that $\mu(\cT_{t}(y^{i, j}_{u1}))=\mu(\cT_{t}(y^{k, \ell}_{w1}))$, $t\ge 0$. 
\end{proof}

Recall the collection of subtrees $\{\cC_{b, i}: i\in \mathrm D_{b}, b\in \Br(\cC)\}$ of $\cC$ indexed by the branch points and their degree ranges. 
Lemma~\ref{lem: cons} implies that for each $\cC_{b, i}$, there is a unique routing point $r^{i, j}_{u1}$ across all $u\in \mathbb U^{\ast}$ and all the pairs $(i, j)$ that is associated to it; we thus denote it as $\sigma_{b, i}$. 
Let $\mathcal W=\{\sigma_{b, i}: i\in \mathrm D_{b}, b\in \Br(\cC)\}$. We are now in a position to prove Propositions~\ref{prop: route} and~\ref{prop: tr-dist}. 

\begin{proof}[Proof of Proposition~\ref{prop: route}]
Let $\mathcal W$ be as defined above. Let $\mathcal R_{\mathcal W}(i, j)=\{\tilde r^{i, j}_{u}: u\in \mathbb U\}$ be the routing for the pair $(i, j)$ associated to $\mathcal W$ as defined in Section~\ref{sec: trace-def}. We show that $\mathcal R_{\mathcal W}(i, j)$ is well-defined and is identical to $\mathcal R(i, j)$. To start off, we have 
\[
\tilde r^{i, j}_{0}=r^{i, j}_{0}=i =  \tilde r^{i, j}_{00}=r^{i, j}_{00} \quad\text{and}\quad \tilde r^{i, j}_{1}=r^{i, j}_{1}= j=  \tilde r^{i, j}_{10}=r^{i, j}_{10}.
\]
Therefore, $\tilde b^{i, j}_{\varnothing}:=\tilde r^{i, j}_{0}\wedge \tilde r^{i, j}_{1}=b^{i, j}_{\varnothing}$. Write $b=b^{i, j}_{\varnothing}$. Then there exists some $k\in \mathrm D_{b}$ so that $\cC_{b, k}$ is the component of $\cC\setminus\{b\}$ containing $i$. By Lemma~\ref{lem: cons} and the definition of $\mathcal W$, we deduce $ \tilde r^{i, j}_{01}=r^{i, j}_{01}$. 
Note that by Lemma \ref{prop: trace-dist'}, $r^{i, j}_{01}$ is independent of $r^{i, j}_{00}$; so the two points are different a.s. Same arguments apply to $ \tilde r^{i, j}_{11}$. We have shown that the elements in $\mathcal R_{\mathcal W}(i, j)$ up to the first generation are well-defined and are identical to those in $\mathcal R(i, j)$. Iterating this argument will allow us to show the same properties for the whole collection. It is now also clear that $\mathcal B_{\mathcal W}(i, j)=\mathcal B(i, j)$, and the second statement in Proposition~\ref{prop: route} then follows from Lemma~\ref{lem: ro-br}. 
\end{proof}

It remains to show Proposition~\ref{prop: tr-dist}. The distributions have been identified in Lemma~\ref{prop: trace-dist'}. It remains to prove the independence. Note this is stronger than the pairwise independence property in Lemma~\ref{prop: trace-dist'}. To prove such a statement, let us first look at its counterpart in $\bp^{n}$-trees. Recall from the paragraph underneath~\eqref{ineq-degree} the sequence of subtrees $(C_{v, i})_{1\le i\le D^{+}_{v}(C)}$ above the vertex $v$. Recall also from there that in each $C_{v, i}$ there is a unique vertex $u_{v, i}$ which was a neighbour of $v$ in $T^{n}$. 
We will need the following result. 
\begin{lem}
\label{lem: p-trace}
The collection $\{u_{v, i}: 1\le i\le D^{+}_{v}(C^{n}), v\in [n]\}$ is independent, and $u_{v, i}\sim \bp^{n}|_{C_{v, i}}$, $1\le i\le D^{+}_{v}(C)$ and $v\in [n]$. 
\end{lem}

\begin{proof}[Proof of Lemma~\ref{lem: p-trace}]
As in Lemma~\ref{lem: p-tree}, this is a consequence of the reconstruction procedure in \cite{BW}. 
\end{proof}

\begin{proof}[Proof of Proposition~\ref{prop: tr-dist}]

Assign a total order $<$ on $\{(i, j)\in \N^{2}: i<j\}$ so that every pair $(i, j)$ appears after finitely many elements in this ordering. 
Note that each branch point $b\in \cC$ can be written as $b=i\wedge j$ for possibly infinitely many pairs of $(i, j)$. Choose the minimal pair with respect to $<$ and call this the minimal writing of $b$. 
Let $(b_{m})_{m\ge 1}$ be the sequence of branch points sorted in the $<$ ordering of their minimal writings. 
Recall that $(V^{n}_{i})_{i\ge 1}$ is a sequence of i.i.d.~vertices of $C^{n}$ with common law $\bp^{n}$, and plays the role of $\N$ in $C^{n}$. Clearly, we can define a similar notion of minimal writing for the branch points of $C^{n}$. Call $(b^{n}_{m})_{1\le m\le \#\Br(C^{n})}$ the sequence of branch points of $C^{n}$, sorted similarly as in $(b_{m})_{m\ge 1}$. 
Under the assumption in~\eqref{p-tree}, the convergence in \eqref{cv: cut-tree} implies that for all $m_{1}<m_{2}<\cdots<m_{k}$, we have
\[
\Big(C^{n}, \alpha_{n}\delta^{n}, \bp^{n}, (b^{n}_{m_{1}},b^{n}_{m_{2}},\dots, b^{n}_{m_{k}} )\Big) \xrightarrow[n\to\infty]{(d)} \Big(\cC, \delta, \nu, (b_{m_{1}}, b_{m_{2}}, \dots, b_{m_{k}})\Big)
\]
with respect to the pointed Gromov--Prokhorov topology. Note that the convergence in \eqref{cv: Cdist} implies that for each $m$ and $i\ge 1$, 
\[
\Big(C^{n}, \alpha_{n}\delta^{n}, \bp^{n}, u_{b^{n}_{m}, i}\Big) \xrightarrow[n\to\infty]{(d)} \Big(\cC, \delta, \nu, \sigma_{b_{m}, i}\Big)
\]
with respect to the pointed Gromov--Prokhorov topology. In fact, we can improve this into a joint convergence: for all $m_{1}<m_{2}<\cdots<m_{k}$ and $j\ge 1$, we have
\[
\Big(C^{n}, \alpha_{n}\delta^{n}, \bp^{n}, \big(u_{b^{n}_{m_{\ell}}, i}\big)_{1\le\ell\le k, 1\le i\le j} \Big) \xrightarrow[n\to\infty]{(d)} \Big(\cC, \delta, \nu, \big(\sigma_{b_{m_{\ell}}, i}\big)_{1\le\ell\le k, 1\le i\le j} \Big)
\]
with respect to the pointed Gromov--Prokhorov topology. The independence of $(\sigma_{b_{m}, i})_{i\in \mathrm D_{b}, b\in \Br(\cT)}$ then follows from Lemma~\ref{lem: p-trace}. 
\end{proof}

{\small
\setlength{\bibsep}{.2em}

}


\begin{thebibliography}{22}
	\providecommand{\natexlab}[1]{#1}
	\providecommand{\url}[1]{\texttt{#1}}
	\expandafter\ifx\csname urlstyle\endcsname\relax
	\providecommand{\doi}[1]{doi: #1}\else
	\providecommand{\doi}{doi: \begingroup \urlstyle{rm}\Url}\fi
	
	\bibitem[Abraham and Delmas(2008)]{AbDe08}
	R.~Abraham and J.-F. Delmas.
	\newblock Fragmentation associated with {L}\'evy processes using snake.
	\newblock \emph{Probab. Theory Related Fields}, 141\penalty0 (1-2):\penalty0
	113--154, 2008.
	
	\bibitem[Abraham and Delmas(2012)]{AbDe12}
	R.~Abraham and J.-F. Delmas.
	\newblock A continuum-tree-valued {M}arkov process.
	\newblock \emph{Ann. Probab.}, 40\penalty0 (3):\penalty0 1167--1211, 2012.
	
	\bibitem[Abraham and Delmas(2013)]{AbDe13}
	R.~Abraham and J.-F. Delmas.
	\newblock The forest associated with the record process on a {L}\'evy tree.
	\newblock \emph{Stochastic Process. Appl.}, 123\penalty0 (9):\penalty0
	3497--3517, 2013.
	
	\bibitem[Abraham et~al.(2010)Abraham, Delmas, and Voisin]{AbDe10}
	R.~Abraham, J.-F. Delmas, and G.~Voisin.
	\newblock Pruning a {L}\'evy continuum random tree.
	\newblock \emph{Electron. J. Probab.}, 15:\penalty0 no. 46, 1429--1473, 2010.
	
	\bibitem[Abraham et~al.(2014)Abraham, Delmas, and Hoscheit]{AbDeHo14}
	R.~Abraham, J.-F. Delmas, and P.~Hoscheit.
	\newblock Exit times for an increasing {L}\'evy tree-valued process.
	\newblock \emph{Probab. Theory Related Fields}, 159\penalty0 (1-2):\penalty0
	357--403, 2014.
	
	\bibitem[Addario-Berry et~al.(2019)Addario-Berry, Dieuleveut, and
	Goldschmidt]{ADG}
	L.~Addario-Berry, D.~Dieuleveut, and C.~Goldschmidt.
	\newblock Inverting the cut-tree transform.
	\newblock \emph{Annales de l'Institut Henri Poincar{\'e} (Probabilit{\'e}s et
		Statistiques)}, 55\penalty0 (3):\penalty0 1349--1376, 2019.
	
	\bibitem[Aldous(1991)]{Aldous1991b}
	D.~Aldous.
	\newblock The continuum random tree. {I}.
	\newblock \emph{The Annals of Probability}, 19:\penalty0 1--28, 1991.
	
	\bibitem[Aldous(1993)]{Aldous1993a}
	D.~Aldous.
	\newblock The continuum random tree {III}.
	\newblock \emph{The Annals of Probability}, 21:\penalty0 248--289, 1993.
	
	\bibitem[Aldous and Pitman(2000)]{Al00}
	D.~Aldous and J.~Pitman.
	\newblock Inhomogeneous continuum random trees and the entrance boundary of the
	additive coalescent.
	\newblock \emph{Probab. Theory Related Fields}, 118\penalty0 (4):\penalty0
	455--482, 2000.
	
	\bibitem[Bertoin(2006)]{Bert_frag}
	J.~Bertoin.
	\newblock \emph{Random fragmentation and coagulation processes}, volume 102 of
	\emph{Cambridge Studies in Advanced Mathematics}.
	\newblock Cambridge University Press, Cambridge, 2006.
	
	\bibitem[Bertoin and Miermont(2013)]{BeMi13}
	J.~Bertoin and G.~Miermont.
	\newblock The cut-tree of large {G}alton-{W}atson trees and the {B}rownian
	{CRT}.
	\newblock \emph{Ann. Appl. Probab.}, 23\penalty0 (4):\penalty0 1469--1493,
	2013.
	
	\bibitem[Bhamidi et~al.(2018)Bhamidi, van~der Hofstad, and Sen]{BHS}
	S.~Bhamidi, R.~van~der Hofstad, and S.~Sen.
	\newblock The multiplicative coalescent, inhomogeneous continuum random trees,
	and new universality classes for critical random graphs.
	\newblock \emph{Probab. Theory Related Fields}, 170\penalty0 (1-2):\penalty0
	387--474, 2018.
	
	\bibitem[Broutin and Wang(2017{\natexlab{a}})]{BW}
	N.~Broutin and M.~Wang.
	\newblock Cutting down $p$-trees and inhomogeneous continuum random trees.
	\newblock \emph{Bernoulli}, 23\penalty0 (4A):\penalty0 2380--2433,
	2017{\natexlab{a}}.
	
	\bibitem[Broutin and Wang(2017{\natexlab{b}})]{BrWa}
	N.~Broutin and M.~Wang.
	\newblock Reversing the cut tree of the brownian continuum random tree.
	\newblock \emph{Electron. J. Probab.}, 22\penalty0 (80):\penalty0 1--23,
	2017{\natexlab{b}}.
	
	\bibitem[Camarri and Pitman(2000)]{Pi00}
	M.~Camarri and J.~Pitman.
	\newblock Limit distributions and random trees derived from the birthday
	problem with unequal probabilities.
	\newblock \emph{Electron. J. Probab.}, 5:\penalty0 no. 2, 18 pp. (electronic),
	2000.
	
	\bibitem[Dieuleveut(2015)]{Di15}
	D.~Dieuleveut.
	\newblock The vertex-cut-tree of {G}alton-{W}atson trees converging to a stable
	tree.
	\newblock \emph{Ann. Appl. Probab.}, 25\penalty0 (4):\penalty0 2215--2262,
	2015.
	
	\bibitem[Duquesne and Le~Gall(2002)]{DuLG02}
	T.~Duquesne and J.-F. Le~Gall.
	\newblock Random trees, {L\'e}vy processes and spatial branching processes.
	\newblock \emph{Ast{\'e}risque}, \penalty0 (281):\penalty0 vi+147, 2002.
	
	\bibitem[Duquesne and Le~Gall(2005)]{DuLeG05}
	T.~Duquesne and J.-F. Le~Gall.
	\newblock Probabilistic and fractal aspects of {L\'e}vy trees.
	\newblock \emph{Probab. Theory Related Fields}, 131\penalty0 (4):\penalty0
	553--603, 2005.
	
	\bibitem[Evans(2008)]{Evans}
	S.~N. Evans.
	\newblock \emph{Probability and real trees}, volume 1920 of \emph{Lecture Notes
		in Mathematics}.
	\newblock Springer, Berlin, 2008.
	
	\bibitem[Greven et~al.(2009)Greven, Pfaffelhuber, and Winter]{GPW09}
	A.~Greven, P.~Pfaffelhuber, and A.~Winter.
	\newblock Convergence in distribution of random metric measure spaces
	({$\Lambda$}-coalescent measure trees).
	\newblock \emph{Probab. Theory Related Fields}, 145\penalty0 (1-2):\penalty0
	285--322, 2009.
	
	\bibitem[Pitman(1999)]{Pi99}
	J.~Pitman.
	\newblock Coalescent random forests.
	\newblock \emph{J. Combin. Theory Ser. A}, 85\penalty0 (2):\penalty0 165--193,
	1999.
	
	\bibitem[Wang(2022+)]{Wa22+}
	M.~Wang.
	\newblock Stable trees as mixings of inhomogeneous continuum random trees.
	\newblock arXiv:2211.07253, 2022.
	
\end{thebibliography}
\end{document}